\documentclass[12pt,a4paper]{article}
\usepackage{amsmath}
\usepackage{amsfonts}
\usepackage{amssymb}
\usepackage[mathcal]{euscript}
\usepackage{amsopn}
\usepackage{amsthm}
\usepackage{amstext}

\usepackage{graphics}

\newcommand {\noqed}
 {\renewcommand {\qed} {}}

\usepackage{smart}
\SectioningParameter{presecnum}{section}=\noexpand\S;
\RestoreLaTeXeqno
\link{equation}{section}

\usepackage{bbm}

\newcommand {\supplus}{\mathop{{\supset}\llap{\raise
0.5pt\hbox{\normalfont\small+}\hskip 0.5pt}}}

\newcommand {\subplus}{\mathop{{\subset}\llap{\raise
0.5pt\hbox{\normalfont\small+}\hskip 0.5pt}}}

\newcommand {\divby}  {\lower 0.15ex \hbox{\,\vdots\,}}
\newcommand {\Aee}    {{\mathbb  A}}
\newcommand {\Cee}    {{\mathbb  C}}
\newcommand {\Eee}    {{\mathbb  E}}
\newcommand {\Fee}    {{\mathbb  F}}
\newcommand {\Gee}    {{\mathbb  G}}
\newcommand {\Hee}    {{\mathbb  H}}
\newcommand {\Kee}    {{\mathbb  K}}
\newcommand {\Nee}    {{\mathbb  N}}
\newcommand {\Oee}    {{\mathbb  O}}
\newcommand {\Pee}    {{\mathbb  P}}
\newcommand {\Qee}    {{\mathbb  Q}}
\newcommand {\Ree}    {{\mathbb  R}}

\newcommand {\Uee}    {{\mathbb  U}}
\newcommand {\Zee}    {{\mathbb  Z}}

\newcommand {\fa}     {{\mathfrak{a}}}
\newcommand {\fA}     {{\mathfrak{A}}}
\newcommand {\fab}    {{\mathfrak{ab}}}
\newcommand {\fag}    {{\mathfrak{ag}}}
\newcommand {\fas}    {{\mathfrak{as}}}
\newcommand {\faut}   {{\mathfrak{aut}}}
\newcommand {\fb}     {{\mathfrak{b}}}
\newcommand {\fc}    {{\mathfrak{c}}}
\newcommand {\fcg}    {{\mathfrak{cg}}}
\newcommand {\fcsvect} {{\mathfrak{csvect}}}
\newcommand {\fcvect}   {{\mathfrak{cvect}}}
\newcommand {\fD}     {{\mathfrak{D}}}   %
\newcommand {\fder}   {{\mathfrak{der}}}   %
\newcommand {\fdiff}  {{\mathfrak{diff}}}
\newcommand {\fdg}    {{\mathfrak{dg}}}
\newcommand {\fd}     {{\mathfrak{d}}}
\newcommand {\fe}     {{\mathfrak{e}}}
\newcommand {\feg}    {{\mathfrak{eg}}}
\newcommand {\ff}     {{\mathfrak{f}}}
\newcommand {\fF}     {{\mathfrak{F}}}
\newcommand {\fg}     {{\mathfrak{g}}}    %
\newcommand {\fgl}    {{\mathfrak{gl}}}  %
\newcommand {\fGL}    {{\mathfrak{GL}}}  %
\newcommand {\fh}     {{\mathfrak{h}}}
\newcommand {\fhei}   {{\mathfrak{hei}}}
\newcommand {\fii}    {{\mathfrak{i}}}    %
\newcommand {\fI}    {{\mathfrak{I}}}    %
\newcommand {\fk}     {{\mathfrak{k}}}

\newcommand {\fl}     {{\mathfrak{l}}}
\newcommand {\fL}     {{\mathfrak{L}}}
\newcommand {\fle}    {{\mathfrak{le}}}
\newcommand {\fM}     {{\mathfrak{M}}}   %
\newcommand {\fm}     {{\mathfrak{m}}}
\newcommand {\fN}     {{\mathfrak{N}}}  %
\newcommand {\fn}     {{\mathfrak{n}}}
\newcommand {\fns}     {{\mathfrak{ns}}}
\newcommand {\fo}     {{\mathfrak{o}}}
\newcommand {\fosp}   {{\mathfrak{osp}}}
\newcommand {\fp}    {{\mathfrak{p}}}   %
\newcommand {\fpe}    {{\mathfrak{pe}}}   %
\newcommand {\fpg}    {{\mathfrak{pg}}}
\newcommand {\fpm}    {{\mathfrak{pm}}}
\newcommand {\fpgl}   {{\mathfrak{pgl}}}
\newcommand {\fpo}    {{\mathfrak{po}}}
\newcommand {\fPo}    {{\mathfrak{Po}}}
\newcommand {\fpsl}   {{\mathfrak{psl}}}
\newcommand {\fpq}    {{\mathfrak{pq}}}
\newcommand {\fpsh}   {{\mathfrak{psh}}}
\newcommand {\fpsq}   {{\mathfrak{psq}}}
\newcommand {\fq}     {{\mathfrak{q}}}
\newcommand {\fr}     {{\mathfrak{r}}}
\newcommand {\fqdiff} {{\mathfrak{qdiff}}}
\newcommand {\fs}     {{\mathfrak{s}}}
\newcommand {\fsa}    {{\mathfrak{sa}}}
\newcommand {\fsb}    {{\mathfrak{sb}}}
\newcommand {\fsg}    {{\mathfrak{sg}}}
\newcommand {\fsh}    {{\mathfrak{sh}}}
\newcommand {\fsi}    {{\mathfrak{si}}}
\newcommand {\fsl}    {{\mathfrak{sl}}}
\newcommand {\fsle}   {{\mathfrak{sle}}}
\newcommand {\fsm}    {{\mathfrak{sm}}}
\newcommand {\fsp}    {{\mathfrak{sp}}}
\newcommand {\fspe}   {{\mathfrak{spe}}}
\newcommand {\fspo}   {{\mathfrak{spo}}}
\newcommand {\fsq}    {{\mathfrak{sq}}}
\newcommand {\fss}    {{\mathfrak{ss}}}
\newcommand {\fst}    {{\mathfrak{st}}}
\newcommand {\fsvect} {{\mathfrak{svect}}}
\newcommand {\fsu}    {{\mathfrak{su}}}
\newcommand {\ft}     {{\mathfrak{t}}}
\newcommand {\fu}     {{\mathfrak{u}}}
\newcommand {\fv}     {{\mathfrak{v}}}     %
\newcommand {\fvect}  {{\mathfrak{vect}}}   %
\newcommand {\fvir}   {{\mathfrak{vir}}}
\newcommand {\fw}     {{\mathfrak{w}}}
\newcommand {\fwitt}  {{\mathfrak{witt}}}
\newcommand {\fx}  {{\mathfrak{x}}}
\newcommand {\fz}     {{\mathfrak{z}}}

\newcommand {\cA}     {{\cal A}}
\newcommand {\cB}     {{\cal B}}
\newcommand {\cC}     {{\cal C}}
\newcommand {\cD}     {{\cal D}}
\newcommand {\cE}     {{\cal E}}
\newcommand {\cF}     {{\cal F}}
\newcommand {\cG}     {{\cal G}}
\newcommand {\cH}     {{\cal H}}
\newcommand {\cI}     {{\cal I}}
\newcommand {\cJ}     {{\cal J}}
\newcommand {\cK}     {{\cal K}}
\newcommand {\cL}     {{\cal L}}
\newcommand {\cM}     {{\cal M}}
\newcommand {\cN}     {{\cal N}}
\newcommand {\cO}     {{\cal O}}
\newcommand {\cP}     {{\cal P}}
\newcommand {\cR}     {{\cal R}}
\newcommand {\cS}     {{\cal S}}
\newcommand {\cT}     {{\cal T}}
\newcommand {\cU}     {{\cal U}}
\newcommand {\cV}     {{\cal V}}
\newcommand {\cW}     {{\cal W}}
\newcommand {\cX}     {{\cal X}}
\newcommand {\cZ}     {{\cal Z}}

\def \opname#1#2%
  {\expandafter\newcommand \csname #1\endcsname {{\mathop{#2}\nolimits}}}

\newcommand{\rmname}[1]
  {\expandafter\newcommand \csname #1\endcsname {{\operatorname{#1}}}}

\newcommand{\rmnameii}[2]
  {\expandafter\newcommand \csname #1\endcsname {{\operatorname{#2}}}}

\rmname{act} \rmname{Ad} \rmname{Add} \rmname{ad} \rmname{Alt}
\rmname{alt} \rmname{Ann} \rmname{antidiag} \rmname{Ber}
\rmname{ber} \rmname{Bil} \rmname{Br} \rmname{card} \rmname{ch}
\rmname{Char} \rmname{cem} \rmname{cj} \rmname{Cliff}
\rmname{cntr} \rmname{codim} \rmname{coind} \rmname{const}
\rmname{col} \rmname{cork} \rmname{cpr} \rmname{diag}
 \rmnameii{Div}{div} \rmname{Def} \rmname{Deg}
\rmname{Der} \rmname{Diff} \rmname{Dim} \rmname{End} \rmname{Even}
\rmname{Ext} \rmname{gr} \rmname{Hom} \rmname{HT}
\rmnameii{Ht}{ht} \rmname{hwt} \rmname{Id} \rmname{id}
\rmname{ind} \rmname{Ind} \rmname{Inf} \rmname{irr} \rmname{Le}
\rmname{Lie} \rmname{lwt} \rmname{mult} \rmname{Mat} \rmname{Mor}
\rmname{nm} \rmname{Ob} \rmname{Odd} \rmname{Osc} \rmname{per}
\rmname{Pic} \rmname{pr} \rmname{pro} \rmname{Prime} \rmname{Proj}
\rmname{prt} \rmname{pt} \rmname{Q} \rmname{qet} \rmname{qtr}
\rmname{rd} \rmname{rk} \rmname{row} \rmname{Res} \rmname{salt}
\rmname{Sch} \rmname{SBr} \rmname{sdim}\rmname{scalar}
\rmname{Ser} \rmname{sign} \rmname{Smbl} \rmname{spin}
\rmname{ssym} \rmname{str} \rmname{st} \rmname{sgn} \rmname{sq}
\rmname{symm} \rmname{supp} \rmname{Supp} \rmname{St}
\rmname{Spec} \rmname{Spm} \rmname{tr} \rmname{vpt} \rmname{Vect}
\rmname{weyl} \rmname{Weyl} \rmname{Witt}

\opname{vvol}  {{v\hspace{-0.1ex}o\hspace{-0.02ex}l\/}}
\opname{pnt}  {\text{\normalfont pt}} \opname{Span} {{Span}}
\opname{slim} {\overline{\lim}} \opname{Vol}
{{V\hspace{-0.55ex}o\hspace{-0.02ex}l\/}} \opname{Par}
{{P\hspace{-0.3ex}a\hspace{-0.05ex}r\/}}

\newcommand {\ev} {{\bar0}}
\newcommand {\od} {{\bar1}}

\newcommand {\tto} {\longrightarrow}
\newcommand {\pder}[1] {{\frac{\partial}{\partial {#1}}}}
\newcommand {\pderf}[2] {{\frac{\partial {#1}}{\partial {#2}}}}

\newcommand {\bcdot}   {\mathbin{\hbox{\raise.4ex\hbox{\bf.}}}} 

\newcommand {\secno} {}
\newcommand {\ssecfont} {\normalfont\bf}

\makeatletter

\theoremstyle{plain}
\@@newtheorem*{Theorem}{\secno Theorem}
\@@newtheorem*{Lemma}{\secno Lemma}
\@@newtheorem*{Proposition}{\secno Proposition}
\@@newtheorem*{Corollary}{\secno Corollary}
\@@newtheorem*{Statement}{\secno Statement}
\@@newtheorem*{Problem}{\secno Problem}
\@@newtheorem*{Question}{\secno Question}
\@@newtheorem*{Conjecture}{\secno Conjecture}

\theoremstyle{definition}
\@@newtheorem*{Example}{\secno Example}
\@@newtheorem*{Examples}{\secno Examples}
\@@newtheorem*{Convention}{\secno Convention}
\@@newtheorem*{Comment}{\secno Comment}

\theoremstyle{remark}

\@@newtheorem*{Remark}{\secno Remark}
\@@newtheorem*{Remarks}{\secno Remarks}
\@@newtheorem*{Exercise}{\secno Exercise}
\@@newtheorem*{Solution}{\secno Solution}
\@@newtheorem*{Hint}{\secno Hint}

\makeatother

\newcommand{\ssec}[1]{\subsection{\boldmath #1.}}

\newcommand {\ssbegin}[1]
  {\refstepcounter{subsection}
  \def \secno {\gdef \secno {}{\ssecfont \thesubsection\hskip 2ex}%
  }%
   \begin{#1}}

\begin{document}
\subjclass{17B50} \keywords {Lie algebras, Lie superalgebras}


\address{MPIMiS, Inselstr. 22, DE-04103 Leipzig, Germany; lebedev@mis.mpg.de}

\title
{Non-degenerate bilinear forms in characteristic $2$, related
contact forms, simple Lie algebras and superalgebras}

\author{Alexei Lebedev
\thanks{I am thankful to the International Max Planck Research School for financial support and most
creative environment and Dimitry Leites for raising the problem
and help.} }
\date{}

\maketitle

\begin{abstract}
Non-degenerate bilinear forms over fields of characteristic $2$,
in particular, non-symmetric ones, are classified with respect to
various equivalences, and the Lie algebras preserving them are
described. Although it is known that there are two series of
distinct finite simple Chevalley groups preserving the
non-degenerate symmetric bilinear forms on the space of even
dimension, the description of simple Lie algebras related to the
ones that preserve these forms is new. The classification of
1-forms is shown to be related to one of the considered
equivalences of bilinear forms. A version of the above results for
superspaces is also given.
\end{abstract}

\section{Introduction}

\ssec{Notations} The ground field $\Kee$ is assumed to be
algebraically closed unless specified; its characteristic is
denoted by $p$; we assume that $p=2$ unless specified;  vector
spaces $V$ are finite dimensional; $n=\dim V$. We often use the
following matrices
\begin{equation}\label{matrices}
\renewcommand{\arraystretch}{1.4}
\begin{array}{l}
\begin{array}{l}
J_{2n}=\begin{pmatrix}0&1_n\\
-1_n&0\end{pmatrix},\quad \Pi_n=\begin{cases} \begin{pmatrix}0&1_k\\
1_k&0\end{pmatrix}&\text{if $n=2k$},\\
\begin{pmatrix}0&0&1_k\\
0&1&0\\
1_k&0&0\end{pmatrix}&\text{if $n=2k+1$},\end{cases}\\
S(n)=\antidiag_n(1, \dots, 1),\quad
Z(2k)=\diag_k(\Pi_2,\dots,\Pi_2).
\end{array}
\end{array}\end{equation}
We call a square matrix {\it zero-diagonal} if it has only zeros
on the main diagonal; let $ZD(n)$ be the space (Lie algebra if
$p=2$) of symmetric zero-diagonal $n\times n$-matrices.

For any Lie (super)algebra $\fg$, let $\fg^{(1)}=[\fg,\fg]$ and
$\fg^{(i+1)}=[\fg^{(i)},\fg^{(i)}]$.

Let $\fo_I(n)$, $\fo_\Pi(n)$ and $\fo_S(n)$ be Lie algebras that
preserve bilinear forms $1_n$, $\Pi_n$ and $S(n)$, respectively.

\ssec{Motivations} Recall that to any bilinear form $B$ on a
given space $V$ one can assign its {\it Gram matrix} by abuse of
notations also denoted by $B=(B_{ij})$: in a fixed basis
$x_1,\dots, x_n$ of $V$ we set
\begin{equation}\label{Gram}
B_{ij}=B(x_i,x_j).
\end{equation}
In what follows, we fix a basis of $V$ and
identify a bilinear form with its matrix. Two bilinear forms $B$
and $C$ on $V$ are said to be {\it equivalent} if there exists an
invertible linear operator $A\in GL(V)$ such that
$B(x,y)=C(Ax,Ay)$ for all $x,y\in V$; in this case,
\begin{equation}\label{eq}
B=ACA^T
\end{equation}
for
the matrices of $B,C$ and $A$ in the same basis.

A bilinear form $B$ on $V$ is said to be {\it symmetric} if $B(v,
w)=B(w, v)$ for any $v, w\in V$; it is {\it skew-symmetric} if
$B(v, v)=0$ for any $v\in V$.

Given a bilinear form $B$, let
$$
L(B)=\{F\in \End\,V\mid B(F x, y)+B(x, Fy)=0\}.
$$
be the Lie algebra that preserves $B$. If $p\neq 2$, some of the
Lie algebras $L(B)$ are simple, for example, the orthogonal Lie
algebras $\fo_B(n)$ that preserve non-degenerate symmetric forms
and symplectic Lie algebras $\fsp_B(n)$ that preserve
non-degenerate skew-symmetric forms.

If $p=2$, either the derived algebras of $L(B)$ for non-degenerate
forms $B$ or their quotient modulo center are simple, so the
canonical expressions of the forms $B$ are needed as a step in
classification of simple Lie algebras in characteristic 2 which is
an open problem, and as a step in a version of this problem for
Lie superalgebras, also open.

The problem of describing preserved bilinear forms has two levels:
we can consider linear transformations (Linear Algebra) and
arbitrary coordinate changes (Differential Geometry). In the
literature, both levels are completely investigated, except for
the case where $p=2$.

More precisely,  for $p=2$, there are obtained rather esoteric
results such as classifications of quadratic forms over skew
fields \cite{ET}, and of analogs of hermitian forms in infinite
dimensional spaces \cite{Gr}, whereas (strangely enough for such a
classically formulated problem) the non-degenerate bilinear forms
over fields were never classified, except for symmetric forms.
Moreover, the fact that the non-split and split forms of the Lie
algebras that preserve the symmetric forms are not always
isomorphic was never mentioned (although known on the Chevalley
group level), cf. the latest papers with reviews of earlier
results \cite{Br, Sh} and \cite{GG}.

Hamelink \cite{H} considered simple Lie algebras over $\Kee$ but
under too restrictive conditions (he considered only Lie algebras
with a nonsingular invariant form) and so missed the fact that
there are two types of orthogonal (or symplectic) simple Lie
algebras.

The bilinear forms over fields of characteristic 2 also naturally
appear in topological problems related to real manifolds, for
example, in singularity theory: as related to \lq\lq symplectic
analogs of Weyl groups" and related bilinear forms over $\Zee/2$,
cf. \cite{I}.

We also consider the superspaces. The Lie superalgebras over
$\Zee/2$ were of huge interest in 1960s in relation with other
applications in topology, see, e.g., \cite{Ha, May}; lately, the
interest comes back \cite{V}.

Let us review the known results and compare them with the new ones
(\S\S 3--8).

\ssec{Known facts: The case $p\neq 2$} Having fixed a basis of
the space on which bilinear or quadratic form is considered, we
identify the form with its Gram matrix; this is understood
throughout. Let me recall the known (both well known and not so
well known) facts. First, recall the elementary Linear Algebra
\cite{Pra}, \cite{L}. Any nondegenerate bilinear form $B$ on a
finite dimensional space $V$ can be represented as the sum $B=S+K$
of a symmetric and a skew-symmetric form. Classics investigated
$B=S+K$ by considering it as a member of the {\it pencil}
$B(\lambda, \mu)=\lambda S+\mu K$, where $\lambda, \mu\in\Kee$,
and studying invariants of $B(\lambda, \mu)$, cf. \cite{Ga}.

Now let the form $B=B(x)$ depend on a parameter $x$ running over a
(super)manifold. Locally, there are obstructions to reducing
non-degenerate 2-form $B(x)$ on a (super)manifold to the canonical
expression. These obstructions are the Riemann tensor if $B(x)$ is
symmetric (metric) and $dB$ if $B(x)$ is a skew (differential)
2-form over $\Cee$ or $\Ree$; for these obstructions expressed in
cohomological terms, see \cite{LPS}. Analogous obstructions to
local reducibility  of bilinear forms to canonical expressions
found here will be classified elsewhere.

\underline{Over $\Cee$}, there is only one class of symmetric
forms and only one class of skew-symmetric forms (\cite{Pra}). For
a canonical form of the matrix of the form $B$, one usually takes
$J_{2n}$ for the skew forms and  $1_n$, $\Pi_n$ or $S(n)$ for
symmetric forms, see \cite{FH}.

In order to have Cartan subalgebra of the orthogonal Lie algebra
on the main diagonal (to have a split form of $\fo_B(n)$), one
should take $B$ of the shape $\Pi_n$ or $S(n)$, not $1_n$. Over
$\Ree$, the Lie algebra might have no split form; for purposes of
representation theory, it is convenient therefore to take its form
most close to the split one.

\underline{Over $\Ree$}, as well as over any ordered field,
Sylvester's theorem states that the signature of the form is the
only invariant (\cite{Pra}).

\underline{Over an algebraically closed field $\Kee$ of
characteristic $p\neq 2$}, Ermolaev  considered
non\-de\-ge\-ne\-rate bilinear forms $B: V\times V\tto \Kee$, and
gave the following description of the Lie algebras $L(B)$ (for
details, see \cite{Er}):

\ssbegin{Statement}[\cite{Er}] The Lie algebra $L(B)$ can not be
represented as a direct sum of ideals (of type $L(B)$) if and only
if all elementary divisors of the matrix $B$ belong to the same
point of the variety $$ P=(\Kee^*/
\{1,-1\})\cup\{0\}\cup\{\infty\}, $$ where $\Kee^*$ is the
multiplicative group of $\Kee$. To each of the three different
types of points in $P$ (elementary divisors corresponding to 0, to
$\infty$ or to a point of $\Kee^*/\{1,-1\}$), a series of Lie
algebras $L(B)$ corresponds, and each of these algebras depends on
a finite system of integer parameters.
\end{Statement}

\ssec{Known facts: The case $p= 2$} 1) With any symmetric
bilinear form $B$ a quadratic form $Q(x):=B(x, x)$ is associated.
The other way round, given a quadratic form $Q$, we define a
symmetric bilinear form, called {\it the polar form} of $Q$, by
setting
$$
B_Q(x, y)=Q(x+y)-Q(x)-Q(y).
$$
As we will see, the correspondence $Q\longleftrightarrow B_Q$ is
not one-to-one and does not embrace non-symmetric forms.

Arf \cite{Arf} has discovered {\it the Arf invariant} --- an
important invariant of non-degenerate quadratic forms in
characteristic 2; for an exposition, see \cite{D}. Two such forms
are equivalent if and only if their Arf invariants are equal.

The Arf invariant, however, can not be used for classification of
symmetric bilinear forms because one symmetric bilinear form can
serve as the polar form for two non-equivalent (and having
different Arf invariants) quadratic forms. Moreover, not every
symmetric bilinear form can be represented as a polar form.

2) Albert \cite{A} classified symmetric bilinear forms over a
field of characteristic 2 and proved that

\medskip
(1) two alternate symmetric forms (he calls a form $B$ on $V$ {\it
alternate} if $B(x,x)=0$ for every $x\in V$; here we call such
forms {\it fully isotropic}) of equal ranks are equivalent;

(2) every non-alternate symmetric form has a matrix which is
equivalent to a diagonal matrix;

(3) if $\Kee$ is such that every element of $\Kee$ has a square
root, then every two non-alternate symmetric forms of equal ranks
are equivalent.

\medskip

3) Albert also gave some results on the classification of
quadratic forms over a field $\Kee$ of characteristic 2
(considered as elements of the quotient space of all bilinear
forms by the space of symmetric alternate forms). In particular,
he showed that if $\Kee$ is algebraically closed, then every
quadratic form is equivalent to exactly one of the forms
\begin{equation}
\label{aforms} x_1x_{r+1}+\dots+x_rx_{2r}\;\text{ or }\;
x_1x_{r+1}+\dots+x_rx_{2r}+x_{2r+1}^2,
\end{equation}
where $2r$
is the rank of the form.

4) Albert also considered {\it semi-definite} bilinear forms,
i.e., symmetric forms, which are equivalent to forms whose matrix
is of the shape
$$
\left(\begin{array}{ll}
1_k&0\\
0 &0\end{array}\right).
$$

For $p=2$, semi-definite forms constitute a linear space. In order
not to have every non-alternate symmetric form semi-definite, one
should take ground field $\Kee$ such that not every element of
$\Kee$ has a square root. For this, $\Kee$ must be neither
algebraically closed nor finite.

5) Skryabin \cite{Sk} considered the case of the space $V$ with a
flag
$$\cF: 0=V_0\subset V_1\subset\dots \subset V_q=V,
$$
and the equivalence of bilinear forms which, in addition to
(\ref{eq}), {\it preserves $\cF$.} He showed that under such
equivalence the class of a (possibly, degenerate) skew-symmetric
bilinear form is determined by parameters
$$
n_{qr}=\dim (V_q\cap V_{r-1}^{\bot})/(V_q\cap
V_r^{\bot}+V_{q-1}\cap V_{r-1}^{\bot})
$$
for $q,r\geq 1$, where orthogonality is taken with respect to the
form. This is true for any characteristic, but if $p=2$, the skew
forms do not differ from zero-diagonal symmetric ones.

\ssec{The structure of the paper} In \S 2 we reproduce Albert's
results on classification of symmetric bilinear forms with respect
to the classical equivalence (\ref{eq}).

In \S 3 we consider other approaches to the classification of
non-symmetric bilinear forms, select the most interesting one
(\lq\lq sociological") and describe the corresponding equivalence
classes.

In \S 4 we classify bilinear forms on superspaces with respect to
the classical and sociological equivalences.

In \S 5 we describe some relations between equivalences of
bilinear forms and 1-forms.

In \S 6 we explicitly describe canonical forms of symmetric
bilinear forms, related simple Lie algebras, and their Cartan
subalgebras.

In \S 7 and \S 8, we give a super versions of \S 6 and \S 3,
respectively.

\ssec{Remarks} 1) In the study of simple Lie algebras over a field
of characteristic $p>0$, one usually takes an algebraically closed
or sometimes finite ground field. Accordingly, these are the cases
where bilinear or quadratic forms are to be considered first.

2) For quadratic forms in characteristic 2, we can also use
Bourbaki's definition: $q$ is quadratic if $q(ax)=a^2q(x)$, and
$B(x,y)=q(x+y)-q(x)-q(y)$ is a bilinear form.

3) For some computations, connected with Lie algebras of linear
transformations, preserving a given bilinear form (e.g.,
computations of Cartan prolongs), it is convenient to choose the
form so that the corresponding Lie algebra has a Cartan subalgebra
as close to the algebra of diagonal matrices as possible. It is
shown in \S 6 that in the case of bilinear forms, equivalent to
$1_n$, over a space of even dimension, we need to take the form in
one of the shapes:
$$
\begin{pmatrix} 1_2&0\\0&S(n-2)\end{pmatrix}\qquad \text{or}\qquad \begin{pmatrix}
1_2&0&0\\0&0&1_{k-1}\\0&1_{k-1}&0\end{pmatrix}
$$

The corresponding Cartan subalgebras consist of matrices of the
following shape:
$$\left(\begin{array}{ll|l} 0&a_0&0\\
a_0&0&0\\
\hline
 0&0 &\diag_n(a_1,\dots,a_{k-1},a_{k-1},\dots,
a_1)\end{array}\right)
$$
or
$$
\left(\begin{array}{ll|l} 0&a_0&0\\
a_0&0&0\\
\hline
 0&0 &\diag_n(a_1,\dots,a_{k-1},a_{1},\dots,
a_{k-1})\end{array}\right).
$$

\section{Symmetric bilinear forms (Linear Algebra)}

\ssbegin{Theorem}\label{SymForm} Let $\Kee$ be a field of
characteristic $2$ such that every element of $\Kee$ has a square
root\footnote{Since $a^2-b^2=(a-b)^2$ if $p=2$, it follows that no
element can have two distinct square roots.}. Let $V$ be a
$n$-dimensional space over $\Kee$.

1) For $n$ odd, there is only one equivalence class of
non-degenerate symmetric bilinear forms on $V$.

2) For $n$ even, there are two equivalence classes of
non-degenerate symmetric bilinear forms, one contains $1_n$ and
the other one contains $S(n)$.
\end{Theorem}

Later we show that, if $n$ is even, a non-degenerate bilinear form
is equivalent to $S(n)$ if and only if its matrix is
zero-diagonal.

Observe that the fact that the bilinear forms are not equivalent
does not imply that the Lie (super)algebras that preserve them are
not isomorphic; therefore the next Lemma is non-trivial.

\begin{Lemma}\label{noniso} 1) The Lie algebras $\fo_I(2k)$ and $\fo_S(2k)$ are not
isomorphic; the same applies to their derived algebras:

2) $\fo_I^{(1)}(2k)\not\simeq \fo_S^{(1)}(2k)$;

3) $\fo_I^{(2)}(2k)\not\simeq \fo_S^{(2)}(2k)$.\end{Lemma}

\ssec{Proof of Theorem \ref{SymForm}} In what follows let
$E^{ij}$, where $1\le i,j \le n$, be a matrix unit, i.e.,
$(E^{i,j})_{kl}:=\delta_{ik}\delta_{jl}$, and
$$
T^{i,j}:=I+E^{i,i}+E^{j,j}+E^{i,j}+E^{j,i}.
$$
 Note that the
$T^{i,j}$ are invertible, and $T^{i,j}=T^{j,i}=(T^{i,j})^T$.

Note also, that any bilinear form $B$ is equivalent to $aB$ for
any $a\in\Kee$ if $a\neq 0$. Indeed, since every element of $\Kee$
has a square root, $aB=(b1_n)B(b1_n)$, where $b\in\Kee$ is such
that $b^2=a$.

Now, let us first prove the following

\ssbegin{Lemma}\label{L1} Let $B$ be a symmetric $n\times n$ matrix, and
$\overline{B}$ be $n'\times n'$ matrix in the upper left corner of
$B$, $n'<n$. Then, if $\overline{B}$ is invertible, $B$ is
equivalent to a matrix of the form
$$
\left(\begin{array}{cc} \overline{B}& 0\\
0&
\widehat{B}\end{array}\right).
$$
\end{Lemma}

\begin{proof}  Let $C$ be $(n-n')\times n'$ matrix in the
lower left corner of $B$, and
$$
M=\left(\begin{array}{ll} 1_{n'}&
0\\C\overline{B}^{-1}&1_{n-n'}\end{array}\right).
$$
The matrix $M$ is invertible because it is lower-triangular and
has no zeros on the diagonal. Direct calculations show that the
matrix $MBM^T$ has the needed form. \end{proof}

\ssbegin{Lemma}\label{L2}  If $B$ and $C$ are $n\times n$ matrices, and
$$
B=\left(\begin{array}{ll} B_1& 0\\0& B_2\end{array}\right);\qquad
C=\left(\begin{array}{ll} C_1& 0\\0& C_2\end{array}\right),
$$
where $B_1$ and $C_1$ are equivalent $n'\times n'$ matrices,
 and $B_2$ and $C_2$ are equivalent $(n-n')\times (n-n')$
matrices, then $B$ and $C$ are equivalent.\end{Lemma}

\begin{proof} If $M_1B_1M_1^T=C_1$ and $M_2B_2M_2^T=C_2$, and
$M=\diag(M_1, M_2)$, then $MBM^T=C$. \end{proof}

Let us fix terminology. A bilinear form $B$ is said to be {\it
fully isotropic} (Albert called them {\it alternate}) if the
corresponding quadratic form $Q(x)=B(x,x)$ is identically equal to
$0$.

\ssbegin{Lemma}\label{ZD} A bilinear form is fully isotropic if and
only if its matrix is zero-diagonal.\end{Lemma}

\begin{proof} Let $e_1,\dots,e_n$ be the basis in which matrix is
taken. Then, if $B$ is fully isotropic, $B_{ii}=B(e_i,e_i)=0$. On
the other hand, if the matrix of $B$ is zero-diagonal, and
$e=\sum_i c_ie_i$, then
$$B(e,e)=\sum_{i,j} B_{ij}c_ic_j=2\sum_{i<j}
B_{ij}c_ic_j=0.\qed
$$
\noqed\end{proof}

Since a fully isotropic form can be equivalent only to a fully
isotropic form, we have

\ssbegin{Corollary}\label{L3} If matrices $A$ and $B$ are symmetric and
equivalent, and $A$ is zero-diagonal, then $B$ is
zero-diagonal.\end{Corollary}

Now, let us prove the following part of Theorem \ref{SymForm}:

\medskip
{\sl If $n=2k$, any non-degenerate symmetric zero-diagonal matrix
$B\in GL(n)$ is equivalent to the matrix $Z(2k)$. }

\medskip

We will prove a more general statement that will be needed later:

\ssbegin{Lemma}\label{ZDForm} If $B$ is a zero-diagonal $n\times n$ matrix (possibly, degenerate), then
$r=\text{rank~} B$ is even, and $B$ is equivalent to the matrix
$$
\tilde{Z}(n,r)=\left(\begin{array}{ll} Z(r)& 0\\0&
0\end{array}\right).
$$
\end{Lemma}

\begin{proof} We will induct on $n$. In the cases $n=1,2$, the statement is
evident.

If $B=0$, the statement follows immediately. Otherwise, there
exist $i,j$ such that $B_{ij}\neq 0$, and $B$ is equivalent to
$$
C=(B_{i,j})^{-1}T^{2,j}T^{1,i}BT^{1,i}T^{2,j},
$$
and $C_{12}=C_{21}=1$, $C_{11}=C_{22}=0$.

Then, by Lemma \ref{L1}, $C$ is equivalent to a matrix $D$ of the
form
$$
\left(\begin{array}{cc|c} 0&1&0\\
1&0&0\\
\hline
 0&0 &D_1\end{array}\right),
 $$
where $D_1$ is a  $(n-2)\times (n-2)$ matrix. Since $D$ is, by
Corollary \ref{L3}, symmetric and zero-diagonal, $D_1$ is also
symmetric and zero-diagonal, and $\text{rank~}D_1=r-2$. Then, by
the induction hypothesis, $r-2$ is even, and $D_1$ is equivalent
to $\tilde{Z}(n-2,r-2)$, and, by Lemma~\ref{L2}, $D$ is equivalent
to $\tilde{Z}(n,r)$. Therefore, $B$ is equivalent to
$\tilde{Z}(n,r)$.\end{proof}

Now let us prove the following:

\medskip
\noindent{\sl For $n$ odd, any non-degenerate symmetric $n\times
n$ matrix is equivalent to $1_n$; \\
for $n$ even, any
non-degenerate symmetric $n\times n$ matrix which is not
zero-diagonal is also equivalent to $1_n$.}

\medskip We will prove this by induction on $n$
(simultaneously for $n$ odd and even). For $n=1$, the statement is
evident.

Now, let  $n$ be even. If $B$ is an invertible symmetric $n\times
n$ matrix, $B_{ii}\neq 0$, then $B$ is equivalent to
$C=(B_{ii})^{-1}T^{1,i}BT^{1,i}$, and $C_{11}=1$. Then, by Lemma
\ref{L1}, $C$ is equivalent to matrix $D$ of the form
\begin{equation}
\label{*}
 \left(\begin{array}{l|l}
1&0\\
\hline 0 &D_1\end{array}\right),\quad \text{where
$D_1\in\fgl(n-1)$.}
\end{equation}
Since $B$ is symmetric and non-degenerate, $D_1$ is also symmetric
and non-degenerate. Then, by induction hypothesis, $D_1$ is
equivalent to $1_{n-1}$, and, by Lemma \ref{L2}, $D$ is equivalent
to $1_n$, and $B$ is also equivalent to $1_n$.

If $B$ is an invertible symmetric $n\times n$ matrix, and $n$ is
odd, then by Lemma \ref{ZDForm}, $B$ has at least one non-zero
element on the diagonal, and, similarly, we can show that $B$ is
equivalent to a matrix $D$ of the form (\ref{*}). Since $B$ is
symmetric and non-degenerate, $D_1$ is also symmetric and
non-degenerate. Then, by induction hypothesis, $D_1$ is equivalent
to either $1_{n-1}$ or $Z(n-1)$, and, by Lemma \ref{L2}, $B$ is
equivalent to either $1_n$ or
$$
\widehat{Z}(n)=\left(\begin{array}{l|l}
1&0\\
\hline 0 &Z(n-1)\end{array}\right).
$$

Let $M$ be a $n\times n$ matrix such that
$$
M_{ij}=\begin{cases}
1&\text{if } i=1\text{~or~} j=1,\\
&\text{or if } j=i,\\
&\text{or if } j=i+1,~i\text{~is odd},\\
&\text{or if } j>i+1;\\
0&\text{if } j=i+1,~i\text{~is even},\\
&\text{or if } 1<j<i.
\end{cases}
$$

Direct calculation shows, that $MM^T=\widehat{Z}(n)$, so
$\widehat{Z}(n)$ is equivalent to $1_n$, and $B$ is equivalent to
$1_n$. 

Now, to finish the proof of the theorem, we need to show that, for
$n$ even, $1_n$ and $Z(n)$ are not equivalent, which follows from
Corollary \ref{L3}. Theorem \ref{SymForm} is proved. \qed

\ssec{Proof of Lemma \ref{noniso}} Let $C(\fg)$ be the center of
the Lie algebra $\fg$. We see that $1_n\in C(\fo_I(n))$, and $\dim
C(\fo_I(n))=1$, because if $A\in \fo_I(n)$, and $A_{ii}\neq
A_{jj}$, then $[A,E^{i,j}+E^{j,i}]_{ij}=A_{ii}+A_{jj}\neq 0$, and
if $A_{ij}\neq 0$ for $i\neq j$, then $[A,E^{i,i}]_{ij}=A_{ij}\neq
0$. Since matrices from $\fo_I(n)^{(1)}$ are zero-diagonal ones,
$\dim (C(\fo_I(n))\cap \fo_I(n)^{(1)})=0$. At the same time,
$$
1_n,\quad \sum_{i=1}^k E^{2i-1,2i},\quad \sum_{i=1}^k
E^{2i,2i-1}\in \fo_S(n); \quad\text{ and } 1_n\in C(\fo_S(n)),
$$
but $[\sum_{i=1}^k E^{2i-1,2i},\quad \sum_{i=1}^k
E^{2i,2i-1}]=1_n$, so $\dim (C(\fo_S(n))\cap \fo_S(n)^{(1)})\neq
0$, which shows that $\fo_I(n)$ and $\fo_S(n)$ are not isomorphic.
Lemma is proved. \qed

\section{Non-symmetric bilinear forms (Linear algebra) }\label{NonSym}

\ssec{Non-symmetric bilinear forms: Discussion} If $p=2$, there is
no canonical way to separate symmetric part of a given bilinear
form from its non-symmetric part, so in this subsection $B$ is
just a non-symmetric form. In this subsection I list several more
or less traditional equivalences before suggesting (in the next
subsection) the one that looks the best.

1) The standard definition (\ref{eq}). This equivalence is too
delicate: there are too many inequivalent forms: the
classification problem looks wild.

2) The idea of classics (see, e.g., \cite{Ga}) was to consider the
following equivalence of non-degenerate bilinear forms regardless
of their symmetry properties. Observe that any bilinear form $B$
on $V$ can be considered as an operator
\begin{equation}
\label{ekJB}
\widetilde{B}:V\tto V^*
 \qquad x\longmapsto B(x,\cdot).
\end{equation}
If $B$ is non-degenerate, then $\widetilde{B}$ is invertible. Two
forms $B$ and $C$ are said to be {\it roughly equivalent}, if the
operators $\widetilde{B}^{-1}\widetilde{B}^*$ and
$\widetilde{C}^{-1}\widetilde{C}^*$ in $V$ are equivalent; here *
denotes the passage to the adjoint operator. This equivalence is,
however, too rough: in the case of symmetric forms it does not
differ fully isotropic and not fully isotropic forms, so {\it all}
symmetric non-degenerate bilinear forms are roughly equivalent,
for both odd- and even-dimensional $V$.

3) Leites suggested to call two bilinear forms $B_1$ and $B_2$
{\it Lie-equivalent} (we  write $B_1\simeq_{L}B_2$) if the Lie
algebras that preserve them are isomorphic. This does reduce the
number of non-equivalent forms but only slightly as compared with
(\ref{eq}) and no general pattern is visible, see the following
Examples for $n=2, 3, 4$. So this equivalence is also, as
(\ref{eq}), too delicate.

\begin{Examples}
\underline{$n=2$, $\Kee=\Zee/2$.} In this case, there exist only
two non-symmetric non-degenerate matrices:
$$
\left(\begin{array}{cc} 1&1\\0&1\end{array}\right)  \qquad
\text{and} \qquad \left(\begin{array}{cc}
1&0\\1&1\end{array}\right),
$$
and they are equivalent.

\underline{$n=2$, $\Kee$ infinite.} In this case, there exist {\it
infinitely many} equivalence classes of non-symmetric
non-degenerate forms. For example,
$\antidiag_2(1,a)\sim\antidiag_2(1,b)$ only if either $a=b$ or
$ab=1$. But all these classes are Lie-equivalent: any
non-symmetric non-degenerate $2\times 2$ matrix is only preserved
by scalar matrices.

\underline{$n=3$, $\Kee=\Zee/2$.} In this case, there exist 3
equivalence classes with the following representatives:
$$
\left(\begin{array}{ccc} 1&1&0\\0&1&0\\0&0&1\end{array}\right),
\qquad \left(\begin{array}{ccc}
1&1&0\\0&1&1\\0&0&1\end{array}\right), \qquad
\left(\begin{array}{ccc} 0&0&1\\0&1&0\\1&1&0\end{array}\right).
$$
Note that the last two matrices are equivalent as forms over an
 extension of $\Zee/2$ with 4 elements. All these matrices are
Lie-equivalent~--- the corresponding Lie algebras are 2-dimensional and, since
they contain $1_3$, commutative.

\underline{$n=3$, $\Kee$ infinite.} In this case, again, there
exist infinitely many equivalence classes.

\begin{Conjecture} All symmetric non-degenerate $3\times 3$
matrices are Lie-equivalent and the corresponding Lie algebras are
$2$-dimensional and commutative. \end{Conjecture}

\underline{$n=4$, $\Kee=\Zee/2$.} In this case, there exist 8
equivalence classes with the following representatives: \footnotesize
$$
B_1=\left(\begin{array}{cccc}
1&1&0&0\\0&1&0&0\\0&0&1&0\\0&0&0&1\end{array}\right);\;\;\;
B_2=\left(\begin{array}{cccc}
1&1&0&0\\0&1&0&0\\0&0&1&1\\0&0&0&1\end{array}\right);\;\;\;
B_3=\left(\begin{array}{cccc}
1&1&0&0\\0&1&1&0\\0&0&1&1\\0&0&0&1\end{array}\right);\;\;\;
B_4=\left(\begin{array}{cccc}
1&1&0&0\\0&1&1&0\\0&0&1&0\\0&0&0&1\end{array}\right);
$$
$$
B_5=\left(\begin{array}{cccc}
0&0&0&1\\0&0&1&0\\0&1&0&0\\1&1&0&0\end{array}\right);\;\;\;
B_6=\left(\begin{array}{cccc}
0&0&0&1\\0&0&1&0\\0&1&1&0\\1&0&1&1\end{array}\right);\;\;\;
B_7=\left(\begin{array}{cccc}
0&0&0&1\\0&0&1&0\\1&0&0&0\\1&1&0&1\end{array}\right);\;\;\;
B_8=\left(\begin{array}{cccc}
0&0&0&1\\0&1&0&0\\0&1&1&0\\1&0&0&0\end{array}\right);
$$
\normalsize

The matrices in the pairs $(B_1, B_4)$, $(B_3, B_7)$, $(B_5, B_6)$
are Lie-equivalent, so there are 5 Lie-equivalence classes.

\underline{$n=4$, $\Kee$ infinite.} Again, there exist infinitely
many equivalence classes. There exist also at least 5
Lie-equivalence classes, described in the previous case.

\end{Examples}

\ssec{A sociological approach to bilinear forms} Instead of
considering non-symmetric forms individually, we can consider the
quotient space $NB(n)$ of the space of all forms modulo the space
of symmetric forms. We will denote the element of this quotient
space with representative $B$, by $\{B\}$. We say that $\{B\}$ and
$\{C\}$ are {\it equivalent}, if there exists an invertible matrix
$M$ such that
$$
\{MBM^T\}=\{C\},\text{ i.e., if $MBM^T-C$ is symmetric}
$$
 (this
definition does not depend on the choice of representatives $B$
and $C$).

Any $\{B\}$ has both degenerate and non-degenerate
representatives: the representative with non-zero elements only
above the diagonal (such representative is unique and
characterizes $\{B\}$) is degenerate, and if we add the unit
matrix to it, we get a non-degenerate representative of $\{B\}$.

Note that $\{B\}$ can be also characterized by the symmetric
zero-diagonal matrix $B+B^T$. The rank of $B+B^T$ is said to be
{\it the rank of $\{B\}$}. According to Lemma \ref{ZDForm}, it is
always even. One can show that it is equal to doubled minimal rank
of representatives of $B$.

\ssbegin{Theorem}\label{ASEq} The forms $\{B\}$ and $\{C\}$ are
equivalent if and only of they have equal ranks.
\end{Theorem}

\begin{proof} Let $\text{rank}\{B\}=\text{rank}\{C\}$, i.e.,
$\text{rank~}(B+B^T)=\text{rank~}(C+C^T)$. Since $B+B^T$ and
$C+C^T$ are zero-diagonal, they are, according to Lemma
\ref{ZDForm}, equivalent, i.e., there exists non-degenerate matrix
$M$ such that $M(B+B^T)M^T=C+C^T$. Then $MBM^T+C=(MBM^T+C)^T$ is
symmetric, and $\{B\}$ is equivalent to $\{C\}$. These arguments
are reversible. \end{proof}

So, we see that $NB(n)$ has $[n/2]+1$ equivalence  classes
consisting of elements with ranks $0,2,\dots, 2[n/2]$. As the
representatives of these classes we can take
$\{\tilde{S}^{n,m}\}$, where $m=0,\dots,[n/2]$ and where
$\tilde{S}^{n,m}$ is $n\times n$ matrix such that
$$
\tilde{S}^{n,m}=\begin{array}{ll} &\begin{array}{ll} n-m&
m\end{array}\\ \begin{array}{c} m\\n-m\end{array}&\left(
\begin{array}{cc} 0&~~S(m)~~\\0&0\end{array}\right)\end{array}
$$

The following definition of the linear transformations preserving
an element of $NB(n)$ seems to be the most natural:
\begin{equation}\label{def}
\text{$X$ preserves $\{B\}$ if $XB+BX^T$ is symmetric.}
\end{equation}

Since
$$
XB+BX^T+(XB+BX^T)^T=X(B+B^T)+(B+B^T)X^T,
$$
$X$ preserves $\{B\}$ if and only if $X$ preserves $B+B^T$. Hence
the transformations preserving $\{B\}$ do form a Lie algebra. One
can check that they also form a Lie algebra if $p\neq 2$, and this
algebra is the Lie algebra of transformations preserving the
skew-symmetric representative of $\{B\}$.

The Lie algebra $\fo_{\tilde{S}^{n,m}}$ is spanned by the matrices
$$
\left(\begin{array}{c|c|c}
A&D&B\\
\hline
0&E&0\\
\hline
 C&F&S(k)AS(k)\end{array}\right),\;\begin{matrix}\text{where $A\in\fgl(k)$; $B, C\in\fgl(k)$ are}\\
\text{such that
$B=S(k)BS(k)$ and $C=S(k)CS(k)$;}\\
\text{$D$ and $F$ are $(n-2k)\times k$ matrices; $E\in\fgl
(n-2k)$.}\end{matrix}
$$
This Lie algebra is isomorphic to the semi-direct sum (the ideal
on the right)
$$
\left(\fo_S(2m)\oplus \fgl(n-2m)\right)\supplus (R_{\fo}\otimes
R_{\fgl}^*),
$$
where $R_{\fo}$ and $R_{\fgl}$ are the spaces of the identity
representations of $\fo_S(2m)$ and $\fgl(n-2m)$, respectively.

If $\text{rank} \{B\}<n$, then the Lie algebra of linear
transformations preserving $\{B\}$ is isomorphic to the Lie
algebra of linear transformations preserving symmetric {\it
degenerate} matrix $B+B^T$. So, it seems natural to call $\{B\}$
{\it non-degenerate} if and only if $\text{rank} \{B\}=n$.

By Theorem \ref{ASEq}, all non-degenerate elements of $NB(n)$
(they only exist if $n$ is even) are equivalent. The Lie algebra
of linear transformations preserving any of these forms is
isomorphic to $\fo_S(2k)$.

\section{Bilinear forms on superspaces (Linear algebra)}

\ssec{Canonical expressions of symmetric bilinear forms on
superspaces and the Lie superalgebras that preserve
them}\label{SCanon} For general background related to Linear
Algebra in superspaces and proofs of the statements of this
subsection, see \cite{LS}.

Speaking of superspaces we denote parity by $\Pi$, and
superdimension by $\sdim$. The operators and bilinear forms are
represented by {\it supermatrices} which we will only consider
here in the {\it standard format}. Recall only that to any
bilinear form $B$ on a given space $V$ one can assign its {\it
Gram matrix} also denoted $B=(B_{ij})$ : in a fixed basis
$x_1,\dots, x_n$ of $V$, we set (compare with (\ref{Gram}); for
$p=2$,  the sign disappears)
\begin{equation}
\label{suGram} B_{ij}=(-1)^{\Pi(B)\Pi(x_i)}B(x_i,x_j).
\end{equation}
In what follows, we fix a basis of $V$ and identify a bilinear
form with its matrix. Two bilinear forms $B$ and $C$ on $V$ are
said to be {\it equivalent} if there exists an invertible {\it
even} linear operator $A\in GL(V)$ such that $B(x,y)=C(Ax,Ay)$ for
all $x,y\in V$; in this case, $B=ACA^T$ for the matrices of $B,C$
and $A$ in the same basis.

Generally, the symmetry of the bilinear forms involves signs which
leads to the notion of supertransposition of the corresponding
Gram supermatrices, but for $p=2$ the supertransposition turns
into transposition.

Recall also that, over superspaces, the parity change sends
symmetric forms to skew-symmetric and the other way round, so if
$p\neq 2$, it suffices to consider only symmetric forms.

In super setting, over $\Cee$, there is only one class of
non-degenerate even symmetric (ortho-symplectic) forms and only
one class of non-degenerate odd symmetric (periplectic) forms in
each superdimension. Over $\Ree$, the invariants of even symmetric
forms are pairs of invariants of the restriction of the form onto
the even and odd subspaces and only one class of odd symmetric
(periplectic) forms.

In order to have Cartan subalgebra of the ortho-symplectic Lie
superalgebra on the main diagonal (to have a split form of
$\fosp(n|2m)$), one should take for the canonical form of $B$ the
expression
$$
\begin{pmatrix}S(n)&0\\0&J_{2m}\end{pmatrix}.
$$

If $p=2$, a given symmetric even form on a superspace can be
represented as a direct sum of two forms on the even subspace and
the odd subspace. For each of these forms, Theorem~\ref{SymForm}
is applicable. Lemma \ref{L2} is also applicable in this case, so
every even symmetric non-degenerate form on a superspace of
dimension $(n_{\overline{0}}|n_{\overline{1}})$ over a field of
characteristic 2 is equivalent to a form of the shape (here:
$i=\bar 0$ or $\bar 1$)
$$
B=\left(\begin{array}{ll}
B_{\overline{0}}&0\\0&B_{\overline{1}}\end{array}\right), \quad
\text{where $B_i=\begin{cases}1_{n_i}&\text{if $n_i$ is
odd;}\\
\text{$1_{n_i}$ or $Z(n_i)$}&\text{if $n_i$ is even.}\end{cases}$}
$$

Before we pass to Lie superalgebras preserving a bilinear form
$B$, let us define the Lie superalgebras for $p=2$. We do not know
any general definition (the definition of squaring given below is
hardly meaningful if the ground field $\Kee$ is finite), but for
an algebraically closed $\Kee$ we are working with it is OK.

A {\it Lie superalgebra} for $p=2$ is a superspace $\fg$ such that
$\fg_\ev$ is a Lie algebra, $\fg_\od$ is an $\fg_\ev$-module (made
into the two-sided one by symmetry) and on $\fg_\od$ a {\it
squaring} (roughly speaking, the halved bracket) is defined
\begin{equation}\label{squaring}
\begin{array}{c}
x\mapsto x^2\quad \text{such that $(ax)^2=a^2x^2$ for any $x\in
\fg_\od$ and $a\in \Kee$, and}\\
(x+y)^2-x^2-y^2\text{~is a bilinear form.}
\end{array}
\end{equation}
Then the bracket of odd elements is defined as following:
$$
{}[x, y]:=(x+y)^2-x^2-y^2
$$
The Jacobi identity for  three odd elements is replaced by the
following relation:
$$
{}[x, x^2]=0 \qquad\text{for any $x\in \fg_\od$}.
$$
The following relation must also be satisfied:
$$
~[x^2,y]=[x,[x,y]]\text{~for any~} x\in\fg_\od, y\in\fg.
$$

The Lie superalgebra preserving $B$ --- by analogy with the
orthosymplectic Lie superalgebras $\fosp$ in characteristic 0 we
call it {\it ortho-orthogonal} and denote $\fo\fo(n_\ev|n_\od)$
--- is spanned by the supermatrices which in the standard format are
of the form
$$
\begin{pmatrix}
A_{\overline{0}}&B_{\overline{0}}C^TB_{\overline{1}}^{-1}\\C&A_{\overline{1}}\end{pmatrix},
\; \begin{matrix}\text{where
$A_{\overline{0}}\in\fo_{B_{\overline{0}}}$,
$A_{\overline{1}}\in\fo_{B_{\overline{1}}}$, and}\\
\text{$C$ is arbitrary $n_{\overline{1}}\times n_{\overline{0}}$
matrix.}\end{matrix}
$$

For an odd symmetric form $B$ on a superspace of dimension
$(n_{\overline{0}}|n_{\overline{1}})$ over a field of
characteristic 2 to be non-degenerate, we need
$n_{\overline{0}}=n_{\overline{1}}=k$, so the matrix of $B$ is of
the shape
$$
\left(\begin{array}{ll} 0&\overline{B}\\
\overline{B}^T&0\end{array}\right),
$$
where $\overline{B}$ is a square invertible matrix. Let us take
$$
M=\left(\begin{array}{ll} 1_k&0\\
0&\overline{B}^{-1}\end{array}\right),
$$
then $B$ is equivalent to
$$
MBM^T=
\left(\begin{array}{ll} 0&1_k\\
1_k&0\end{array}\right).
$$

This form is preserved by linear transformations with
supermatrices in the standard format of the shape
\begin{equation}
\label{pe} \left(\begin{array}{ll} A&C\\D&A^T\end{array}\right),
\quad \text{where $A\in\fgl(k)$, $C$ and $D$ are symmetric
$k\times k$ matrices}. 
\end{equation}

As over $\Cee$ or $\Ree$, the Lie superalgebra $\fpe(n)$ of
supermatrices $(\ref{pe})$ (recall that $p=2$) will be referred to
as {\it periplectic}.

\ssec{Non-symmetric forms on superspaces}\label{non-sym_super} If
a non-symmetric form on a superspace is even, it can be again
represented as a direct sum of two bilinear forms: one on the even
subspace, and the other one on the odd subspace. These two forms
can be independently transformed to canonical forms, see
\ref{NonSym}.

The situation with odd non-symmetric forms is more interesting.
Such a form can be non-degenerate only on a space of
superdimension $(k|k)$. In the standard format, the supermatrix of
such a form has the shape
$$
B=\left(\begin{array}{ll} 0&A\\C&0\end{array}\right),
$$
where $A$ and $C$ are invertible matrices. Let $M$ be an
invertible matrix such that $L=MC(A^T)^{-1}M^{-1}$ is the Jordan
normal form of $C(A^T)^{-1}$. Then $B$ is equivalent to
$$
\left(\begin{array}{ll}(M^T)^{-1}A^{-1}&0\\0&M\end{array}\right)
\left(\begin{array}{ll} 0&A\\C&0\end{array}\right)
\left(\begin{array}{ll}((M^T)^{-1}A^{-1})^T&0\\0&M^T\end{array}\right)=
\left(\begin{array}{ll} 0&1_k\\L&0\end{array}\right).
$$

This expression (with $L$ in the Jordan normal form) can be
considered as a canonical form of a non-degenerate odd bilinear
form.

\begin{Statement}  Two non-degenerate forms are equivalent if and only if they have equal
canonical forms. \end{Statement}

\section{Relation with 1-forms (Differential geometry)}

\ssec{Notations} For $p>0$, there are two types of analogs of
polynomial algebra: the infinite dimensional ones and finite
dimensional ones. The {\it divided power algebra} in
indeterminates $ x_1, \dots, x_m$ is the algebra of polynomials in
these indeterminates, so, as space, it is
$$
\cO(m)=\Span\{x_1^{(r_1)}\dots x_m^{(r_m)}\mid r_1,\dots, r_m\geq
0 \}
$$
with the following multiplication:
$$
(x_1^{(r_1)}\dots x_m^{(r_m)})\cdot(x_1^{(s_1)}\dots
x_m^{(s_m)})=\prod\limits_{i=1}^m \left(\begin{array}{c}
r_i+s_i\\r_i \end{array} \right) x_1^{(r_1+s_1)}\dots
x_m^{(r_m+s_m)}.
$$

For a  {\it shearing parameter} $\underline{N}=(n_1,\dots,n_m)$,
set
$$
\cO(m,\underline{N})=\Span\{x_1^{(r_1)}\dots x_m^{(r_m)}\mid 0\leq
r_i<p^{n_i}, i=1,\dots,m \},
$$
where $p^\infty=\infty$. If $n_i<\infty$ for all $i$, then
$\dim\cO(m,\underline{N})<\infty$.

\ssec{Matrices and 1-forms} Let $B$ and $B'$ be the matrices of
bilinear forms on an $n$-dimensional space $V$ over a field $\Kee$
of characteristic 2. Let $x_0,x_1,\dots,x_n$ be indeterminates as
in sec. 5.1; set
$$
\text{$\deg x_0=2$, $\deg x_1=\dots=\deg x_n=1$.}
$$

We say that $B$ and $B'$ are {\it 1-form-equivalent}  if there
exists a degree preserving transformation, i.e., a set of
independent variables $x'_0,x'_1,\dots,x'_n$ such that
\begin{equation}\label{eqdegs}
\deg x'_0=2, \quad \deg
x'_1=\dots=\deg x'_n=1,
\end{equation}
which are polynomials in $x_0,x_1,\dots,x_n$ in divided powers
with shearing parameter
\begin{equation}\label{shear}
\underline{N}=(N_0,\dots,N_n)\;\text{ such that $N_i>1$ for
every $i$ from $1$ to $n$,}
\end{equation}
and such that
\begin{equation}\label{1FEq}
dx_0+\sum\limits_{i,j=1}^n
B_{ij}x_idx_j=dx'_0+\sum\limits_{i,j=1}^n B'_{ij}x'_idx'_j.
\end{equation}

\ssbegin{Lemma}\label{T1FEq} $B$ and $B'$ are 1-form-equivalent if and only if
$\{B\}$ and $\{B'\}$ are equivalent.\end{Lemma}

\begin{proof} By $(\ref{eqdegs})$,
we have
\begin{equation}\label{XTrans}
x'_0=cx_0+\sum\limits_{i=1}^n A_{ii}x_i^{(2)}+\sum\limits_{1\leq
i<j\leq n} A_{ij}x_ix_j;\qquad x'_i=\sum\limits_{j=1}^n M_{ij}x_j,
\end{equation}
where $c\neq 0$ and $M$ is an invertible matrix. Thanks to
(\ref{1FEq}), comparing coefficients of $dx_0$ in the left- and
right-hand sides, we get $c=1$. Let $A$ be a symmetric $n\times n$
matrix with elements $A_{ij}$ for $i\leq j$ as in (\ref{XTrans}).
Then
$$
dx'_0+\sum\limits_{i,j=1}^nB'_{ij}x'_idx'_j=dx_0+\sum\limits_{i,j=1}^n
A_{ij}x_idx_j+\sum\limits_{i,j,k,l=1}^n M_{ki}B'_{kl}M_{lj}x_i
dx_j,
$$
i.e., $B=M^TB'M+A$, so $\{B\}$ and $\{B'\}$ are equivalent. Since
our arguments are invertible, the theorem is proved.\end{proof}

\ssec{The case of odd indeterminates} Let us modify the
definition of 1-form-equivalence for the super case where
$x_1,\dots,x_n$ are all {\it odd}. In this case, we can only use
divided powers with $\underline{N}=(N_0,1,\dots,1)$.

We say that $B$ and $B'$ are {\it 1-superform-equivalent}  if
there exists a set of indeterminates $x'_0,x'_1,\dots,x'_n$, which
are polynomials in $x_0,x_1,\dots,x_n$, such that
\begin{equation}\label{2feq}
\text{$\Pi(x'_0)=\bar 0$, $\Pi(x'_1)=\dots=\Pi(x'_n)=\bar
1$,\qquad $\deg x'_0=2$, $\deg x'_1=\dots=\deg x'_n=1$}
\end{equation}
and
\begin{equation}\label{1SFEq}
dx_0+\sum\limits_{i,j=1}^n
B_{ij}x_idx_j=dx'_0+\sum\limits_{i,j=1}^n B'_{ij}x'_idx'_j.
\end{equation}

Now, recall that the 1-form $\alpha$ on a superdomain $M$ is said
to be {\it contact} if it singles out a nonintegrable distribution
in the tangent bundle $TM$ and $d\alpha$ is non-degenerate on the
fibers of this distribution; for details, see \cite{GL} and
\cite{LPS}.

\ssbegin{Lemma}\label{T1SFEq} The matrices $B$ and $B'$ are 1-superform-equivalent if and only if
exist an invertible matrix $M$ and a symmetric zero-diagonal
matrix $A$ such that
\begin{equation}\label{1feq} B=MB'M^T+A.
\end{equation}\end{Lemma}

\begin{proof} It is analogous to the proof of Lemma
\ref{T1FEq}.\end{proof}

Albert \cite{A} considered the equivalence (\ref{1feq}) as an
equivalence of (matrices of) quadratic forms. In particular, he
proved the following

\ssbegin{Statement}\label{Albert} If $\Kee$ is algebraically closed, every matrix $B$ is
equivalent in the sense $(\ref{1feq})$ to exactly one of the matrices
$$
Y(n,r)=\left(\begin{array}{lll}
0_r&1_r&0\\
0_r&0_r&0\\
0&0&0_{n-2r}\end{array}\right)\qquad \text{or}\qquad
\widetilde{Y}(n,r)=\left(\begin{array}{llll}
0_r&1_r&0&0\\
0_r&0_r&0&0\\
0&0&1&0 \\
0&0&0&0_{n-2r-1}\end{array}\right),
$$
where $2r=\text{rank} (B+B^T)$. The corresponding quadratic form
is non-degenerate if and only if either $n=2r$, or $n=2r+1$ and
the matrix is equivalent to $\widetilde{Y}(n,r)$.\end{Statement}

If, in 1-form-equivalence, we consider divided powers with
shearing parameter $\underline{N}=(N_0,1,\dots,1)$, it is the same
as to consider 1-superform-equivalence.

\ssbegin{Lemma} Let $x_0,\dots,x_n$ be indeterminates,
$\Pi(x_0)=\bar 0$, $\Pi(x_1)=\dots=\Pi(x_n)=\Pi$. Then the
$1$-form on the $(n+1)$-dimensional superspace
\begin{equation}
\label{Lin1F}
 \alpha=dx_0+\sum\limits_{i,j=1}^n B_{ij}x_idx_j
\end{equation}
is contact if and only if one of the following conditions holds:

1) $\Pi=\bar 0$, and $\{B\}$ is non-degenerate, i.e.,
$n=\text{rank}(B+B^T)$ (this rank is always even);

2) $\Pi=\bar 1$, and the quadratic form corresponding to $B$ is
non-degenerate, i.e., either $n=\text{rank}(B+B^T)\;\text{ or
}\;B\;\text{is not zero-diagonal and}\; n=\text{rank}(B+B^T)+1.$
\end{Lemma}

\begin{proof} From Theorem \ref{ASEq}, Statement \ref{Albert} (\cite{A}), and
Lemmas \ref{T1FEq} and \ref{T1SFEq} we know that, if $\Pi=\bar 0$,
every symmetric bilinear form is 1-form-equivalent to one of the
forms $\widetilde{S}(n,r)$, $n\geq 2r$, and, if $\Pi=\bar 1$,
every symmetric bilinear form is 1-superform-equivalent to one of
the forms $Y(n,r)$, where $n\geq 2r$, or $\widetilde{Y}(n,r)$,
where $n\geq 2r+1$. Direct calculations show that if $\Pi=\bar 0$,
the 1-form (\ref{Lin1F}) corresponding to $\widetilde{S}(n,r)$ is
contact if and only if $n=2r$; if $\Pi=\bar 1$, then the 1-form
(\ref{Lin1F}) corresponding to $Y(n,r)$ is contact if and only if
$n=2r$ and the 1-form, corresponding to $\widetilde{Y}(n,r)$ is
contact if and only if $n=2r+1$. Since, by definition, two 1-forms
that correspond to 1-(super)form-equivalent bilinear forms can be
transformed into each other by a change of coordinates, we are
done.
\end{proof}

From this, we get the following:

\begin{Theorem} The following are the canonical
expressions of the odd contact forms:
\begin{equation}
\label{canon} \alpha=dx_0+\sum\limits_{i=1}^k
x_idx_{k+i}\begin{cases}&\text{for $n=2k$ and $x_1,\dots,x_n$ all even or all odd;}\\
+x_ndx_n&\text{for $n=2k+1$ and $x_1,\dots x_n$ odd.}\end{cases}
\end{equation}
\end{Theorem}

\begin{Remarks} 1) If $n>1$ and $x_1,\dots x_n$ are odd, the 1-form
$\alpha=dx_0+\sum\limits_{i=1}^n x_idx_i$ is not contact since
(recall that $p=2$)
$$
\alpha=d\left(x_0+\sum\limits_{i<j} x_ix_j\right)+\left(\sum_{i=1}^n
x_i\right)d\left(\sum_{i=1}^n x_i\right).
$$

2) Let $p=2$. Since there are two types of orthogonal Lie algebras
if $n$ is even, and orthogonal algebras coincide, in a sense, with
symplectic ones, it seems natural to expect that there are also
two types of the Lie algebras of hamiltonian vector fields
(preserving $I$ and $S$, respectively). Iyer investigated this
question; for the answer, see \cite{Iy}.

Are there two types of contact Lie algebras corresponding to these
cases? The (somewhat unexpected) answer is NO:

The classes of 1-(super)form-equivalence of bilinear forms which
correspond to contact forms have nothing to do with classes of
classical equivalence of symmetric bilinear forms. The 1-forms,
corresponding to symmetric bilinear forms are exact if $x_1,\dots
x_n$ are even, and are of rank $\leq 2$ if $x_1,\dots x_n$ are
odd.

Recall the the contact Lie superalgebra consists of the vector
fields $D$ that preserve the contact structure (nonintegrable
distribution given by a contact form $\alpha$ hereafter in the
form (\ref{canon})). Such fields satisfy
$$
L_D(\alpha)=F_D\alpha\;\text{ for some $F_D\in\Cee [t, p, q,
\theta]$.}
$$
For any $f\in\Cee [t, p, q, \theta]$, we set\index{$K_f$, contact
vector field} \index{$H_f$, Hamiltonian vector field} (the signs
here are important only for $p\neq 2$):
\begin{equation}
\label{2.3.1} K_f=(1-E)(f)\pder{t}-H_f + \pderf{f}{t} E,
\end{equation}
where $E=\sum\limits_i y_i \pder{y_{i}}$ (here the $y_{i}$ are all
the coordinates except $t$) is the {\it Euler operator}, and $H_f$
is the hamiltonian field with Hamiltonian $f$ that preserves
$d\alpha_1$:
\begin{equation}
\label{2.3.2} H_f=\sum\limits_{i\leq n}\left(\pderf{f}{p_i}
\pder{q_i}-\pderf{f}{q_i} \pder{p_i}\right )
-(-1)^{p(f)}\left(\sum\limits_{j\leq m}\pderf{ f}{\theta_j}
\pder{\theta_j}\right ) .
\end{equation}

If one tries to build a contact algebra $\fg$ by means of a
non-degenerate symmetric bilinear form $B$ on the space $V$ by
setting (like it is done in characteristic $0$) $\fg$ to be the
generalized Cartan prolongation $(\fg_-,\fg_0)_*$ (for the precise
definition, see \cite{Shch}), where the non-positive terms of
$\fg$ are (here $K_f$ is the contact vector field with the
generating function $f$; for an exact formula, see, e.g.,
\cite{GLS}):
$$
\renewcommand{\arraystretch}{1.4}
\begin{array}{l}
\fg_i=\begin{cases}0&\text{ if }i\leq
-3;\\
\Kee\cdot K_1&\text{ if } i=-2\\
V&\text{ if } i=-1\\
\fo_B(V)\oplus \Kee K_t &\text{ if }i=0\end{cases}\end{array}
$$
and where the multiplication is given by the formulas
$$
\renewcommand{\arraystretch}{1.4}
\begin{array}{l}
{}[X,Y]=B(X,Y)K_1 \text{~for any~} X,Y\in\fg_{-1};\\
\fo_B(V) \text{~acts on $V$ via the standard action};\\
{}[\fg_0,\fg_{-2}]=0;\\
K_t \text{~acts as $\id$ on $\fg_{-1}$},\\
{}[K_t,\fo_B(V)]=0,
\end{array}
$$
then the form $B$ must be zero-diagonal one (because
$0=[X,X]=B(X,X)K_1$ for $X\in\fg_{-1}$).

One can also try to construct a Lie superalgebra in a similar way
by setting $\Pi(\fg_{-1})=\od$ and
\begin{equation}
\label{g_-1}
 X^2=B(X,X)K_1 \text{~for any~} X\in\fg_{-1}.
\end{equation}

Let us realize this algebra by vector fields on a superspace of
superdimension $(1|n)$ with basis $x_0,\dots,x_n$ such that
$$
\Pi(x_0)=\ev;\qquad \Pi(x_i)=\od \text{~for~}1\leq i\leq n.
$$

If $e_1,\dots,e_n$ is a basis of $V$ and we set (here
$\partial_i=\frac{\partial}{\partial x_i}$ for $i=0,\dots,n$):
$$
K_1=\partial_0;\qquad e_i=\partial_i+\sum\limits_{j=1}^n A_{ij}
x_j\partial_0\text{~for~} i=1,\dots,n,
$$
then, to satisfy relations (\ref{g_-1}), we need the following
(here the Gram matrix $B$ is taken in the basis $e_1,\dots,e_n$):
$$
A_{ii}=B_{ii}\text{~for~}1\leq i\leq n;
$$
$$
A_{ij}+A_{ji}=B_{ij}+B_{ji}\text{~for~}1\leq i<j\leq n
$$
i.e., $A\in\{B\}$, where the equivalence class is taken with
respect to zero-diagonal symmetric matrices.

These vector fields preserve the $1$-form
$$
\alpha=dx_0+\sum\limits_{i,j=1}^n A_{ij}x_idx_j.
$$

So, to get a contact Lie superalgebra in this way, one needs $B$
to be non-symmetric with non-degenerate class $\{B\}$.

3) Lin \cite{LinK} considered an $n$-parameter family of simple
Lie algebras for $p=2$ preserving in dimension $2n+1$ the
distribution given by the contact form
$$
\alpha=dt+\sum\limits_{i=1}^n
\left((1-a_{i})p_idq_i+a_{i}q_idp_i\right),\;\text{where
$a_{i}\in\Kee$}.
$$
Obviously, the linear change
\begin{equation}
\label{Linchange} t'=t+\sum a_{i}p_iq_i\quad\text{and identical on
other indeterminates}
\end{equation}
reduces $\alpha$ to the canonical form $dt+\sum\limits_{i=1}^n
p_idq_i$. So the parameters $a_i$ can be eliminated. Although Lin
mentioned the change (\ref{Linchange}) on p.~21 of \cite{LinK},
its consequence was not formulated and, seven years after, Brown
\cite{Br} reproduced Lin's misleading $n$-parameter description of
$\fk(2n+1)$.
\end{Remarks}

\ssec{The case of indeterminates of different parities}

\subsubsection{The case of an odd 1-form.} Let
$$
\text{$\Pi(x_0)=\Pi(x_1)=\dots=\Pi(x_{n_{\overline{0}}})=\bar
0$,\qquad  $\Pi(x_{n_{\overline{0}}+1})=\dots=\Pi(x_n)=\bar 1$.}
$$
This corresponds to the following equivalence (we call it
1-superform-equivalence again) of {\it even} bilinear forms on a
superspace $V$ of superdimension
$(n_{\overline{0}}|n_{\overline{1}})$, where
$n_{\overline{1}}=n-n_{\overline{0}}$: two such forms $B$ and $B'$
are said to be 1{\it-superform-equivalent} if, for their
supermatrices, we have (\ref{1feq}), where $M\in
GL(n_{\overline{0}}|n_{\overline{1}})$ and $A$ is a symmetric even
supermatrix such that the restriction of the bilinear form
corresponding to it onto the odd subspace $V_\od$ is fully
isotropic. This means that, in the standard format of
supermatrices,
$$
B=\left(\begin{array}{ll}B_{\overline{0}}&0\\0&B_{\overline{1}}\end{array}\right)\qquad
\text{and}\qquad
B'=\left(\begin{array}{ll}B'_{\overline{0}}&0\\0&B'_{\overline{1}}\end{array}\right)
$$
are 1-superform-equivalent if and only if (1) $B_{\overline{0}}$
and $B'_{\overline{0}}$ are 1-form-equivalent, and (2)
$B_{\overline{1}}$ and $B'_{\overline{1}}$ are
1-superform-equivalent. (This also follows from the fact that the
1-form
$$
dx_0+\sum_{i,j=1}^{n_{\overline{0}}} A_{ij}x_idx_j+
\sum_{i,j=1}^{n_{\overline{1}}}
B_{ij}x_{n_{\overline{0}}+i}dx_{n_{\overline{0}}+j}
$$ is contact if
and only if the forms $dx_0+\sum_{i,j=1}^{n_{\overline{0}}}
A_{ij}x_idx_j$ and $dx_0+\sum_{i,j=1}^{n_{\overline{1}}}
B_{ij}x_{n_{\overline{0}}+i}dx_{n_{\overline{0}}+j}$ are contact
on the superspaces of superdimension $(n_{\overline{0}}+1|0)$ and
$(1|n_{\overline{1}})$, respectively.) Then, from (\ref{canon}) we
get the following

\begin{Theorem} The following are the canonical expressions for an odd contact form on a superspace:
$$
dt+\sum\limits_{i=1}^k p_idq_{i}+\sum\limits_{j=1}^l
\xi_id\eta_{i}\begin{cases}&\text{for $n_{\overline{0}}=2k$ and $n_{\overline{1}}=2l$,}\\
+\theta d\theta&\text{for $n_{\overline{0}}=2k$ and
$n_{\overline{1}}=2l+1$},
\end{cases}
$$
where $t=x_0$; $p_i=x_i$, $q_i=x_{k+i}$ for $1\leq i\leq k$;
$\xi_i=x_{n_{\overline{0}}+i}$, $\eta_i=x_{n_{\overline{0}}+l+i}$
for $1\leq i\leq l$; $\theta=x_n$ for $n_{\overline{1}}=2l+1$.
\end{Theorem}

\subsubsection{The case of an even 1-form.} Let $\Pi(x_0)=\bar 1$.
This corresponds to the following equivalence of {\it odd}
bilinear forms on a superspace $V$ of superdimension
$(n_{\overline{0}}|n_{\overline{1}})$: two such forms $B$ and $B'$
are said to be 1-{\it superform-equivalent} if for their
(super)matrices we have (\ref{1feq}), where $M\in
GL(n_{\overline{0}}|n_{\overline{1}})$ and $A$ is a symmetric odd
supermatrix. Then, since
$$
\renewcommand{\arraystretch}{1.4}
\begin{array}{l}
\left(\begin{array}{l|l}
1_{n_{\overline{0}}}&0\\
\hline 0&M\end{array}\right)\left(B+\left(\begin{array}{l|l}
0&C\\
\hline C^T&0\end{array}\right)\right)\left(\begin{array}{l|l}
1_{n_{\overline{0}}}&0\\
\hline 0&M^T\end{array}\right)=\left(\begin{array}{l|l}
0&0\\
\hline
X(D+C^T)&0\end{array}\right)\\
\text{for}\;B=\left(\begin{array}{l|l}
0&C\\
\hline D&0\end{array}\right),\end{array}
$$
any such $B$ is equivalent to a form with a supermatrix of the
shape (the indices above and to the left of the supermatrix are
the sizes of the blocks)
$$
\begin{array}{ll}
&\begin{array}{ll|l}
r&n_{\overline{0}}-r~&n_{\overline{1}}\end{array}\\
\begin{array}{l} n_{\overline{0}}\\ \hline r\\
n_{\overline{1}}-r\end{array}& \left(\begin{array}{ll|l}
0&~~~0~~~&0\\
\hline
 1_r&0&0\\
 0&0&0
 \end{array}\right)\end{array},
$$
where $r=\text{rank}(D+C^T)$. The corresponding form is contact if
and only if $r=n_{\overline{0}}=n_{\overline{1}}$. Hence, we get
the following somewhat unexpected result:

\begin{Theorem} The following expressions for the canonical form of an even
contact (pericontact) 1-form on a superspace of dimension
$(k|k+1)$ are equivalent:
$$
d\tau+\sum\limits_{i=1}^k \xi_i dq_i, \quad \text{or }\;
d\tau+\sum\limits_{i=1}^k q_id\xi_i, \quad \text{or }\;
d\tau+\sum\limits_{i=1}^l \xi_i dq_i+\sum\limits_{i=l+1}^k
q_id\xi_i,
$$
where $\tau=x_0$, $\xi_i=x_{k+i}$, $q_i=x_i$ for $1\leq i\leq k$.
\end{Theorem}

\section{Canonical expressions of symmetric bilinear forms. Related
simple Lie algebras} If we want to have a canonical expression of
a non-degenerate bilinear form $B$ such that the intersection of
the Cartan subalgebra of $\fo_B^{(1)}(n)$ or $\fo_B^{(2)}(n)$ with
the space of diagonal matrices were of maximal possible dimension,
we should take the following canonical forms of $B$. Each of the
following subsections \ref{ss61}, \ref{ss62}, \ref{ss63} contains
two most convenient expressions of an equivalence class of
bilinear forms.

\ssec{$n=2k+1$}\label{ss61}

\subsubsection{}\label{ss611} If $B=S(2k+1)$, then $\fo_B(n)$ consists of the
matrices, symmetric with respect to the side diagonal;  it is
convenient to express them in the block form
$$
\begin{pmatrix}
A&X&C\\Y^TS(K)&z&X^TS(k)\\D&Y&S(k)A^TS(k)\end{pmatrix}
\;\begin{matrix}\text{where $A\in\fgl(k)$, $C$ and $D$ are
symmetric with respect to}\\
\text{the side diagonal, $X, Y\in\Kee^k$ are column-vectors,
$z\in\Kee$.}\end{matrix}
$$
The Lie algebra $\fo_B^{(1)}(n)$ consists of the elements of
$\fo_B(n)$, which have only zeros on the side diagonal; the Cartan
subalgebra of $\fo_B^{(1)}(n)$ of maximal dimension is spanned by
the matrices
$$
\diag_n(a_1,\dots,a_{k},0,a_{k},\dots, a_1).
$$

\subsubsection{}\label{ss612}  If $B=\Pi_{2k+1}$, then $\fo_B(n)$ is spanned by the
matrices
$$
\begin{pmatrix}
A&X&C\\Y^T&z&X^T\\D&Y&A^T\end{pmatrix} \;\begin{matrix}\text{where
$A\in\fgl(k)$, $C$ and $D$ are
symmetric,}\\
\text{$X, Y\in\Kee^k$ are column-vectors, $z\in\Kee$.}\end{matrix}
$$
The Lie algebra $\fo_B^{(1)}(n)$ consists of the elements of
$\fo_B$ such that $C$ and $D$ are zero-diagonal, $z=0$; the Cartan
subalgebra of $\fo_B^{(1)}(n)$ of maximal dimension is spanned by
the matrices
$$
\diag_n(a_1,\dots,a_{k},0, a_1, \dots, a_{k}).
$$

\ssec{$n=2k$ and $B$ equivalent to $S(2k)$}\label{ss62}

\subsubsection{}\label{ss6121}  If $B=S(2k)$, then $\fo_B(n)$ consists of the matrices,
symmetric with respect to the side diagonal;  it is convenient to
express them in the block form
$$
\begin{pmatrix}
A&C\\D&S(k)A^TS(k)\end{pmatrix} \;\begin{matrix}\text{where
$A\in\fgl(k)$, $C$ and $D$ are
symmetric}\\
\text{with respect to the side diagonal.}\end{matrix}
$$
The Cartan subalgebra of the related simple Lie algebra (it is
described later) is spanned by the matrices
$$
\diag_n(a_1,\dots,a_{k},a_{k},\dots, a_1)~\text{such
that}~a_1+\dots+a_k=0.
$$

\subsubsection{}\label{ss622}  If $B=\Pi_{2k}$, then $\fo_B(n)$ is spanned by the matrices
\begin{equation}
\label{**}
\begin{pmatrix} A&C\\D&A^T\end{pmatrix},
\;\text{where $A\in\fgl(k)$, $C$ and $D$ are symmetric.}
\end{equation}

Observe that these matrices can be represented as $\Pi(2k)U$ or
$V\Pi(2k)$, where $U$ and $V$ are symmetric.

The Cartan subalgebra of the related simple Lie algebra (it is
described later) is spanned by the matrices
$$
\diag_n(a_1,\dots,a_{k},a_{1},\dots, a_k)~\text{such
that}~a_1+\dots+a_k=0.
$$

\ssec{$n=2k$ and $B$ equivalent to $1_n$}\label{ss63} We get the
greatest dimension of the intersection of the Cartan subalgebra
with the space of diagonal matrices if the matrix of $B$ is of any
of the following shapes:

\subsubsection{}\label{ss631} If $B=\begin{pmatrix} 1_2&0\\0&S(n-2)\end{pmatrix}$,
then $\fo_B$ is spanned by the matrices
$$
\begin{pmatrix}
A&C\\S(n-2)B^T&D\end{pmatrix} \;\begin{matrix}\text{where
$A\in\fgl(2)$ is symmetric,
$C$ is any $2\times (n-2)$ matrix,}\\
\text{$D\in\fgl(n-2)$ is symmetric with respect to the side
diagonal.}\end{matrix}
$$
The Lie algebra $\fo_B^{(1)}(n)$ consists of the elements of
$\fo_B(n)$ such that $A$ is zero-diagonal, $D$ has only zeros on
the side diagonal; the Cartan subalgebra of $\fo_B^{(1)}(n)$ of
greatest dimension is spanned by the matrices
$$\left(\begin{array}{ll|l} 0&a_0&0\\
a_0&0&0\\
\hline
 0&0 &\diag_n(a_1,\dots,a_{k-1},a_{k-1},\dots,
a_1)\end{array}\right).
$$

\subsubsection{}\label{ss632} If $B=\begin{pmatrix}
1_2&0&0\\0&0&1_{k-1}\\0&1_{k-1}&0\end{pmatrix}$, then $\fo_B(n)$
is spanned by the matrices
$$\begin{pmatrix}
X&Y&Z\\Z^T&A&C\\Y^T&D&A^T\end{pmatrix} \;\begin{matrix}\text{where
$X\in\fgl(2)$ is symmetric, $Y$ and $Z$ are of size $2\times
(k-1)$
,}\\
\text{ $A\in\fgl(k-1)$, \ $C, D\in \fgl(k-1)$ are
symmetric}.\end{matrix}
$$
The Lie algebra $\fo_B^{(1)}(n)$ consists of the elements of
$\fo_B(n)$ such that $X$, $C$ and $D$ are zero-diagonal; the
Cartan subalgebra of $\fo_B^{(1)}(n)$  greatest dimension is
spanned by the matrices
$$
\left(\begin{array}{ll|l} 0&a_0&0\\
a_0&0&0\\
\hline
 0&0 &\diag_n(a_1,\dots,a_{k-1},a_{1},\dots,
a_{k-1})\end{array}\right).
$$

\ssec{The derived algebras of $\fo_I(n)$}\label{fo_I-der}  Direct
calculation shows that
$$\fo_I^{(1)}(n)=\begin{cases}0&\text{if $n=1$}\\
\{\lambda S(2)\mid \lambda\in\Kee\}&\text{if
$n=2$};\end{cases}\qquad \fo_I^{(2)}(n)=0 \text{~if~} n\leq 2.
$$

\begin{Lemma}\label{o_I} If $n>2$, then

i) $\fo_I^{(1)}(n)= ZD(n)$;

ii) $\fo_I^{(2)}(n)=\fo_I^{(1)}(n)$.\end{Lemma}

\begin{proof} First, let us show that $\fo_I^{(1)}(n)\subset
ZD(n)$. Indeed, if $A,A'\in\fo_I(n)$, then
$$
[A,A']_{ii}=\sum\limits_j A_{ij}A'_{ji}-A'_{ij}A_{ji}=0
$$
since $A,A'$ are symmetric. So, matrices from $\fo_I^{(1)}(n)$ are
zero-diagonal.

Let $F^{ij}=E^{ij}+E^{ji}$, where $1\leq i,j\leq n$, $i\neq j$.
These matrices are symmetric, so they are all in $ZD(n)$. Let us
show that they also are in $\fo_I^{(1)}(n)$. Since, for $1\leq
i<j\leq n$, the matrices $F^{ij}$  form a basis of $ZD(n)$, it
follows that $ZD(n)\subset\fo_I^{(1)}(n)$, and it proves (i).

Direct calculation shows that if $1\leq k\leq n$, $k\neq i,j$,
then
\begin{equation}
\label{Fcomm} [F^{ik}, F^{kj}]=F^{ij},
\end{equation}
so $F^{ij}\in\fo_I^{(1)}(n)$.

Moreover, once we have shown that $F^{ij}\in\fo_I^{(1)}(n)$, this
computation also proves that $F^{ij}\in\fo_I^{(2)}(n)$. Since
$\fo_I^{(2)}(n)\subset\fo_I^{(1)}(n)=ZD(n)$, it also proves
(ii).\end{proof}

\begin{Lemma}\label{sim_o_I} If $n>2$, then $\fo_I^{(1)}(n)$ is
simple. \end{Lemma}

\begin{proof} Let $I\subset\fo_I^{(1)}(n)$ be an ideal, and $x\in
I$ an element, such that its decomposition with respect to
$\{F^{ij}\}$ contains $F^{ab}$ with non-zero coefficient for some
$a,b$. Let us note that
\begin{equation}
\label{FFF} [F^{ij}, [F^{ij}, F^{kl}]]=\begin{cases}  F^{kl}&
\text{if~} |\{i,j\}\cap \{k,l\}|=1\\ 0&
\text{otherwise.}\end{cases}
\end{equation}

Let us define an operator $P_{F^{ab}}:\fo_I^{(1)}(n)\to
\fo_I^{(1)}(n)$ as follows:
\begin{equation}
\label{o_I_P} P_{F^{ab}}=\begin{cases}(\ad ~F^{bc})^2(\ad
~F^{ac})^2\text{~~for ~}c\neq a,b,~~1\leq c\leq 3&\text{if~}
n=3 ;\\
\prod\limits_{1\leq c\leq n,~ c\neq a,b}(\ad ~F^{ac})^2
&\text{if~} n>3.\end{cases}
\end{equation}

Then, from (\ref{FFF}),
\begin{equation}
\label{o_I_P_act} P_{F^{ab}}
F^{cd}=\begin{cases}F^{cd}&\text{~if~}
F^{cd}=F^{ab};\\0&\text{~otherwise.}\end{cases}
\end{equation}

So, $[F^{ac},[F^{ac},[F^{bc},[F^{bc},x]]]]$ is proportional (with
non-zero coefficient) to $F^{ab}$, and \\$F^{ab}\in I$. Then, from
(\ref{Fcomm}), $F^{ib}, F^{ij}\in I$ for all $i,j$, $1\leq i,j\leq
n$, $i\neq j$, and $I=\fo_I^{(1)}(n)$. \end{proof}

\ssec{The derived Lie algebras of $\fo_{\Pi}(2n)$}\label{fo_P-der}
Direct computations show that:
\medskip

$\fo_{\Pi}^{(1)}(2)=\{\lambda\cdot 1_2\mid \lambda\in\Kee\}$;

$\fo_{\Pi}^{(2)}(2)=0$;

$\fo_{\Pi}^{(1)}(4)=\{\text{matrices of the shape (\ref{**}) such
that $B, C\in ZD(2)$}\}$;

$\fo_{\Pi}^{(2)}(4)=\{\text{matrices of $\fo_{\Pi}^{(1)}(4)$ such
that $\tr A=0$}\}$;

$\fo_{\Pi}^{(3)}(4)=\{\lambda\cdot 1_4\mid \lambda\in\Kee\}$;

$\fo_{\Pi}^{(4)}(4)=0$.

\begin{Lemma}\label{o_p} If $n\geq 3$, then

i) $\fo_{\Pi}^{(1)}(2n)=\{\text{matrices of the shape (\ref{**})
such that $B, C\in ZD(n)$}\}$;

ii) $\fo_{\Pi}^{(2)}(2n)=\{\text{matrices of the shape (\ref{**})
such that $B, C\in ZD(n)$, and $\tr~ A=0$}\}$;

iii) $\fo_{\Pi}^{(3)}(2n)=\fo_{\Pi}^{(2)}(2n)$. \end{Lemma}

\begin{proof} Let $M^1$ and $M^2$ denote
{\it conjectural} $\fo_{\Pi}^{(1)}(2n)$ and $\fo_{\Pi}^{(2)}(2n)$,
respectively, as described in Lemma. First, let us prove that
$\fo_{\Pi}^{(1)}(2n)\subset M^1$ and $\fo_{\Pi}^{(2)}(2n)\subset
M^2$. Let
$$
L=\begin{pmatrix} A&B\\C&A^T\end{pmatrix}, \qquad
L'=\begin{pmatrix} A'&B'\\C'&A'^T\end{pmatrix}\in\fo(2n),\qquad
\text{and~~} L''=[L,L']=\begin{pmatrix}
A''&B''\\C''&A''^T\end{pmatrix}.
$$

Then, for any $i\in\overline{1,n}$, we have
$$
B''_{ii}=\sum\limits_{j=1}^n
(A_{ij}B'_{ji}+B_{ij}A'_{ij}-A'_{ij}B_{ji}-B'_{ij}A_{ij})=0
$$
since $B$,$B'$ are symmetric. Analogically, $C''_{ii}=0$, so
$L''\in M^1$. Hence, $\fo_{\Pi}^{(1)}(2n)\subset M^1$.

Now, if $L,L'\in \fo_{\Pi}^{(1)}(2n)$, then
$$
\tr A''=\sum\limits_{i=1}^n A''_{ii}=\sum\limits_{i,j=1}^n
(A_{ij}A'_{ji}+B_{ij}C'_{ji}-A'_{ij}A_{ji}-B'_{ij}C_{ji})=0
$$
since $B,B',C,C'$ are symmetric and zero-diagonal. So, since
$L''\in\fo_{\Pi}^{(2)}(2n)\subset\fo_{\Pi}^{(1)}(2n)\subset M^1$,
it follows that $L''\in M^2$, and $\fo_{\Pi}^{(2)}(2n)\subset
M^2$.

Let us introduce the following notations for matrices from
$\fo_{\Pi}(2n)$:

\medskip
$F_1^{ij}$, where $1\leq i,j\leq n$, $i\neq j$, such that\qquad
$A=C=0$, $B=E^{ij}+E^{ji}$;

$F_2^{ij}$, where $1\leq i,j\leq n$, $i\neq j$, such that\qquad
$A=B=0$, $C=E^{ij}+E^{ji}$;

$G^{ij}$, where $1\leq i,j\leq n$, $i\neq j$, such that\qquad
$B=C=0$, $A=E^{ij}$;

$H^{ij}$, where $1\leq i,j\leq n$, $i\neq j$, such that\qquad
$B=C=0$, $A=E^{ii}+E^{jj}$;

$K_0$ such that\qquad $B=C=0$, $A=E^{11}$;

$K_1$ such that\qquad $A=C=0$, $B=E^{11}$;

$K_2$ such that\qquad $A=B=0$, $C=E^{11}$;
\medskip

\noindent Observe that $F_1^{ij}$, $F_2^{ij}$, $G^{ij}$, and
$H^{ij}$ span $M^2$; whereas $M^2$ and $K_0$ span $M^1$.

Direct computations give the following relations:
\begin{equation}
\label{rels}
\begin{array}{l} \text{if $k\neq i,j$, then $[H^{ik},
F_1^{ij}]=F_1^{ij}$,\qquad $[H^{ik},
F_2^{ij}]=F_2^{ij}$,\qquad $[H^{ik}, G^{ij}]=G^{ij}$;}\\
{}[F_1^{ij},F_2^{ij}]=H^{ij};\\
{}[K_1,K_2]=K_0.\end{array} \end{equation}

\noindent Since $F_1^{ij},F_2^{ij},G^{ij},H^{ij},K_1,K_2\in
\fo_{\Pi}(2n)$, it follows that
$F_1^{ij},F_2^{ij},G^{ij},H^{ij},K_0\in \fo^{(1)}_{\Pi}(2n)$.
Hence, $M^1\subset\fo^{(1)}_{\Pi}(2n)$, and
$\fo^{(1)}_{\Pi}(2n)=M^1$. Relations (\ref{rels}) imply that
$M^2\subset [M^2,M^2]$, so $M^2\subset
[M^1,M^1]=\fo^{(2)}_{\Pi}(2n)$, and $\fo^{(2)}_{\Pi}(2n)=M^2$.
Also, $M^2\subset [M^2,M^2]=\fo^{(3)}_{\Pi}(2n)$, so
$\fo^{(3)}_{\Pi}(2n)=M^2$. The lemma is proven.\end{proof}

\begin{Lemma}\label{o_p_sim} If $n\geq 3$, then

i) if $n$ is odd, then $\fo_{\Pi}^{(2)}(2n)$ is simple;

ii) if $n$ is even, then the only non-trivial ideal of
$\fo_{\Pi}^{(2)}(2n)$ is the center\\ $Z=\{\lambda\cdot 1_{2n}\mid
\lambda\in\Kee\}$ (thus, $\fo_{\Pi}^{(2)}(2n)/Z$ is
simple).\end{Lemma}

\begin{proof} We use the notations of the previous Lemma. It
follows from the relations
$$
[F_1^{ij},F_2^{ij}]=H^{ij};
$$
$$
[H^{ij},X^{kl}]=\begin{cases} X^{kl}& \text{if~}
|\{i,j\}\cap\{k,l\}|=1\\ 0& \text{otherwise}\end{cases}\qquad
\text{for~} X^{kl}=F_1^{kl},F_2^{kl},G^{kl};
$$
that if an ideal $I$ of $\fo_{\Pi}^{(2)}(2n)$ contains any of the
elements $F_1^{ij}, F_2^{ij}$, then $I=\fo_{\Pi}^{(2)}(2n)$.

Let $1\leq i,j,k\leq n$, $i\neq j\neq k\neq i$. Direct computation
shows that the operators
$$
P_{F_1^{ij}}=\ad F_1^{jk}\ad F_1^{ij}\ad F_2^{ij}\ad F_2^{jk};
\qquad P_{F_2^{ij}}=\ad F_2^{jk}\ad F_2^{ij}\ad F_1^{ij}\ad
F_1^{jk}
$$
on $\fo_{\Pi}^{(2)}(2n)$ act as follows: for $X$ equal to one of
the elements $F_1^{lm}, F_2^{lm}, H^{lm}, G^{lm}$,
$$
P_{F_1^{ij}}X=\begin{cases}X&\text{if~}X=F_1^{ij}\\0&\text{otherwise}\end{cases};\qquad\qquad
P_{F_2^{ij}}X=\begin{cases}X&\text{if~}X=F_2^{ij}\\0&\text{otherwise}\end{cases}.
$$

It follows from these two facts that any element of a non-trivial
ideal $I$ of $\fo_{\Pi}^{(2)}(2n)$ must not contain $F_1^{ij},
F_2^{ij}$ in its decomposition with respect to the basis of
$F_1^{ij}, F_2^{ij}, H^{ij}, G^{ij}$ --- i.e., it must have the
shape
$$
\begin{pmatrix} A&0\\0&A^T\end{pmatrix}.
$$

Then, for $I$ to be an ideal, $A$ must satisfy the following
condition:
$$
AB+BA^T=0\text{~for all~} B\in ZD(n).
$$

If $A$ contains a non-zero non-diagonal entry $A_{ij}$, then
$$
(A(E^{jk}+E^{kj})+(E^{jk}+E^{kj})A^T)_{ik}=A_{ij}\neq 0
$$
for $k\neq i,j$; if $A$ contains two non-equal diagonal entries
$A_{ii}$ and $A_{jj}$, then
$$
(A(E^{ij}+E^{ji})+(E^{ij}+E^{ji})A^T)_{ij}=A_{ii}-A_{jj}\neq 0
$$

So, $A$ must be proportional to $1_n$, and
$1_{2n}\in\fo_{\Pi}^{(2)}(2n)$ if and only if $n$ is even.
\end{proof}

\ssec{On simplicity of the derived algebras of $\fo_B(n)$} Let
$\fo_B(n)$ be either $\fo_I(n)$ or $\fo_{\Pi}(n)$ or $\fo_{S}(n)$,
the Lie algebras of linear transformations preserving bilinear
forms $1_n$ or $\Pi_n$ or $S(n)$, respectively.

\begin{Lemma}\label{simple}
$\fo_B^{(1)}(n)$ is a simple Lie algebra if $n\equiv 1\pmod 2$ and
$n>1$;

$\fo_B^{(2)}(n)/\mathfrak{center}$ is a simple Lie algebra if
$n\equiv 0\pmod 2$ and $n>4$.
\end{Lemma}

\begin{proof} Follows from the fact that in the $\Zee$-grading
corresponding to the blocks of the matrix forms 6.1-6.3, the
$\fg_0$-modules $\fg_{\pm 1}$ are irreducible, generate $\fg_\pm$,
and $[\fg_1,\fg_{-1}]=\fg_0$.

\end{proof}

\section{Canonical expressions of symmetric bilinear superforms. Related
 Lie superalgebras.} In this section we consider Lie
superalgebras of linear transformations preserving bilinear forms
on a superspace of superdimension $(n_\ev|n_\od)$ and their
derived superalgebras. Since in the case where $n_\ev=0$ or
$n_\od=0$ these superalgebras are entirely even and do not differ
from the corresponding Lie algebras, we do not consider this case.

As it was said in sec.~\ref{SCanon}, every even symmetric
non-degenerate form on a superspace of superdimension
$(n_{\overline{0}}|n_{\overline{1}})$ over a field of
characteristic 2 is equivalent to a form of the shape (here:
$i=\bar 0$ or $\bar 1$)
\begin{equation}
\label{canon_sup_ev}
 B=\left(\begin{array}{ll}
B_{\overline{0}}&0\\0&B_{\overline{1}}\end{array}\right), \quad
\text{where $B_i=\begin{cases}1_{n_i}&\text{if $n_i$ is
odd;}\\
\text{$1_{n_i}$ or $\Pi(n_i)$}&\text{if $n_i$ is
even.}\end{cases}$}
\end{equation}
The Lie superalgebra $\fo\fo_B(n_\ev|n_\od)$ preserving $B$ is
spanned by the supermatrices which in the standard format are of
the shape
\begin{equation}
\label{oo_b}
\begin{pmatrix}
A_{\overline{0}}&B_{\overline{0}}C^TB_{\overline{1}}^{-1}\\C&A_{\overline{1}}\end{pmatrix},
\; \begin{matrix}\text{where
$A_{\overline{0}}\in\fo_{B_{\overline{0}}}(n_\ev)$,
$A_{\overline{1}}\in\fo_{B_{\overline{1}}}(n_\od)$, and}\\
\text{$C$ is arbitrary $n_{\overline{1}}\times n_{\overline{0}}$
matrix.}\end{matrix}
\end{equation}

In what follows we use the fact that in the case of matrices and
supermatrices behave identically with respect to multiplication
and Lie (super)bracket
--- i.e., if two square supermatrices of the same format are
given, then the entries of their product or Lie (super)bracket do
not depend on this format.

\ssec{The derived Lie superalgebras of $\fo\fo_{II}(n_{\ev}|
n_{\od})$} Let $B$ be of the shape (\ref{canon_sup_ev}) such that
$B_i=1_{n_i}$. We will denote Lie superalgebra preserving this
form as $\fo\fo_{II}(n_{\ev}| n_{\od})$; this superalgebra
consists of symmetric supermatrices.

Direct calculation shows that
$$
\fo\fo_{II}^{(i)}(1|1)=\begin{cases}
\left\{\begin{pmatrix}a&b\\b&a\end{pmatrix}\mid
a,b\in\Kee\right \} &\text{ if } i=1,\\
\{a\cdot 1_{1|1}\mid a\in\Kee\}&\text{ if } i=2,\\
0&\text{ if }i\geq 3.
\end{cases}
$$

\begin{Lemma} If $n=n_\ev+n_\od\geq 3$, then

i) $\fo\fo_{II}^{(1)}(n_\ev|n_\od)$ consists of symmetric
supermatrices of (super)trace $0$;

ii)
$\fo\fo_{II}^{(2)}(n_\ev|n_\od)=\fo\fo_{II}^{(1)}(n_\ev|n_\od)$.\end{Lemma}

\begin{proof} It was shown in the proof of Lemma \ref{o_I}
that a (super)bracket of any two symmetric matrices is
zero-diagonal, so to prove that supermatrices from
$\fo\fo_{II}^{(1)}(n_\ev|n_\od)$ have trace $0$, we only need to
prove this for the squares of odd symmetric supermatrices. If $L$
is an odd matrix of the shape (\ref{oo_b}), then
$$
\tr L^2=\sum\limits_{i=1}^{n_\ev}(\sum\limits_{j=1}^{n_\od}
C_{ij})^2+
\sum\limits_{j=1}^{n_\od}(\sum\limits_{i=1}^{n_\ev}C_{ij})^2=
2\sum\limits_{i=1}^{n_\ev}\sum\limits_{j=1}^{n_\od}C_{ij}^2=0.
$$

Now let us introduce the following notations for matrices from
$\fo\fo_{II}(n_\ev|n_\od)$: $$
\renewcommand{\arraystretch}{1.4}
\begin{array}{l}
F^{ij}=E^{ij}+E^{ji}\quad\text{ for $1\leq i,j\leq n$, $i\neq
j$;}\\
H^{ij}=E^{ii}+E^{jj}\quad\text{for $1\leq i\leq n_\ev$, $1\leq
j\leq n_\od$.} \end{array}
$$

These matrices span the space of symmetric matrices with trace 0.
As it was shown in the proof of Lemma \ref{o_I}, if $n\geq 3$,
then the matrices $F^{ij}$ generate themselves. Moreover, if
$1\leq i\leq n_\ev$, $1\leq j\leq n_\od$ (so that $F^{ij}$ is
odd), then $(F^{ij})^2=H^{ij}$. So all the $F^{ij}$, $H^{ij}$ lie
in $\fo\fo_{II}^{(1)}(n_\ev|n_\od)$ for any $i$. \end{proof}

\begin{Lemma} If $n=n_\ev+n_\od\geq 3$, then

i) if $n$ is odd, then $\fo\fo_{II}^{(1)}(n_\ev|n_\od)$ is simple;

ii) if $n$ is even, then the only non-trivial ideal of
$\fo\fo_{II}^{(1)}(n_\ev|n_\od)$ is the center\\ $C=\{\lambda\cdot
1_n\mid \lambda\in\Kee\}$ (thus,
$\fo\fo_{II}^{(1)}(n_\ev|n_\od)/C$ is simple).\end{Lemma}

\begin{proof} Let us define operators $P_{F^{ab}}$ as in
(\ref{o_I_P}). Then, due to (\ref{o_I_P_act}) and the fact that
$$
P_{F^{ab}}H^{kl}=0,
$$
we can show in the same way as in the proof of Lemma \ref{sim_o_I}
that if an ideal of $\fo\fo_{II}^{(1)}(n_\ev|n_\od)$ contains a
non-diagonal matrix, it contains all the $F^{ij}$. Since all the
$H^{ij}$ are squares of odd $F^{ij}$, such an ideal is trivial.

So, any non-trivial ideal of $\fo\fo_{II}^{(1)}(n_\ev|n_\od)$ is
diagonal. For a diagonal matrix $X$,
$$
[X,F^{ij}]=(X_{jj}-X_{ii})F^{ij},
$$
so all the elements of a non-trivial ideal must be proportional to
$1_n$, and $1_n\in\fo\fo_{II}^{(1)}(n_\ev|n_\od)$ if and only if
$n$ is even.\end{proof}

\ssec{The derived Lie superalgebras of $\fo\fo_{I\Pi}(n_{\ev}|
n_{\od})$} Now let us consider the case where $n_\od$ is even and
$B$ is of the shape (\ref{canon_sup_ev}) such that
$B_\ev=1_{n_\ev}$, $B_\od=\Pi(n_\od)$. (The case where $n_\ev$ is
even, $B_\ev=\Pi(n_\ev)$, $B_\od=1_{n_{\od}}$ is analogous to this
one, so we will not consider it.) We will denote Lie superalgebra
preserving this form by $\fo\fo_{I\Pi}(n_{\ev}| n_{\od})$; this
superalgebra consists of supermatrices of the following shape:
\begin{equation}
\label{oo_IS}
\begin{pmatrix}
A_{\overline{0}}&C^T\Pi(n_\od)\\C&\Pi(n_\od)A_{\overline{1}}\end{pmatrix},
\;
\begin{matrix}\text{where
$A_{\overline{0}}$, $A_\od$ are symmetric},\\
\text{$C$ is an arbitrary $n_{\overline{1}}\times
n_{\overline{0}}$ matrix.}\end{matrix}
\end{equation}

\begin{Lemma} i) $\fo\fo_{I\Pi}^{(1)}(n_{\ev}|n_{\od})$ consists of the
matrices of the shape (\ref{oo_IS}) such that $A_\ev$ is
zero-diagonal;

ii)
$\fo\fo_{I\Pi}^{(2)}(n_{\ev}|n_{\od})=\fo\fo_{I\Pi}^{(1)}(n_{\ev}|n_{\od})$.\end{Lemma}

\begin{proof} Set $k_\od=n_\od/2$; let $M$ be the {\it
conjectural} space of $\fo\fo_{I\Pi}^{(1)}(n_{\ev}|n_{\od})$ as it
is described in the Lemma. First, let us prove that
$\fo\fo_{I\Pi}^{(1)}(n_{\ev}|n_{\od})\subset M$. If
$$
\begin{array}{l}
L=\begin{pmatrix}
A_{\overline{0}}&C^T\Pi(n_\od)\\C&A_{\overline{1}}\end{pmatrix},
L'=\begin{pmatrix}
A'_{\overline{0}}&C'^T\Pi(n_\od)\\C'&A'_{\overline{1}}\end{pmatrix}\in
\fo\fo_{I\Pi}(n_{\ev}|n_{\od}),\text{~and~}\\
L''=[L,L']=\begin{pmatrix}
A''_{\overline{0}}&C''^T\Pi(n_\od)\\C''&A''_{\overline{1}}\end{pmatrix},
\end{array}
$$
then
$$
\begin{array}{l}
(A''_\ev)_{ii}=([A_\ev,A'_\ev]+C^T\Pi(n_\od)C'-C'^T\Pi(n_\od)C)_{ii}=
\\ \sum\limits_{j=1}^{n_\ev}
((A_\ev)_{ij}(A'_\ev)_{ji}-(A'_\ev)_{ij}(A_\ev)_{ji})+\sum\limits_{j=1}^{k_\od}
(C_{ji}C'_{j+k_\od,i}+C_{j+k_\od,i}C'_{ji}-C'_{ji}C_{j+k_\od,i}+C'_{j+k_\od,i}C_{ji})=0
\end{array}
$$
since $A_\ev,A'_\ev$ are symmetric. Now, if $L$ is odd (i.e.,
$A_\ev=0$, $A_\od=0$), then
$$
L=\begin{pmatrix} C^T\Pi(n_\od)C&0\\0&CC^T\Pi(n_\od)\end{pmatrix},
$$
and
$$
(C^T\Pi(n_\od)C)_{ii}=\sum\limits_{j=1}^{k_\od}
(C_{ji}C_{j+k_\od,i}+C_{j+k_\od,i}C_{ji})=0.
$$

Let us introduce the following notations for matrices from
$\fo\fo_{I\Pi}(n_{\ev}|n_{\od})$:
\medskip

$F^{ij}$, where $1\leq i,j\leq n_\ev$ and $i\neq j$, such that
$C=0$, $A_\od=0$, $A_\ev=E^{ij}+E^{ji}$;

$G^{ij}$, where $1\leq i\leq n_\od$, $1\leq j\leq n_\ev$, such
that $A_\ev=0$, $A_\od=0$, $C=E^{ij}$;

$H^i$, where $1\leq i\leq n_\od$, such that $A_\ev=0$, $C=0$,
$A_\od=\begin{cases} E^{i,i+k_\od}& \text{~if~} 1\leq i\leq
k_\od;\\E^{i,i-k_\od}& \text{~if~} k_\od+1\leq i\leq
n_\od;\end{cases}$

$I_{\od}$, such that $A_\ev=0$, $C=0$, $A_\od=1_{n_\od}$.
\medskip

\noindent Direct calculations give the following relations:
$$
\begin{array}{l}
{}[G^{i1},G^{k_\od+1,j}]=F^{ij};\\
{}[I_\od, G^{ij}]=G^{ij};\\
{}(G^{i1})^2=H^i,
\end{array}
$$
so we get $F^{ij},G^{ij},H^i\in
\fo\fo_{I\Pi}^{(1)}(n_{\ev}|n_{\od})$.

Let us also denote by $K$ the subalgebra of
$\fo\fo_{I\Pi}(n_{\ev}|n_{\od})$, consisting of all matrices of
the form (\ref{oo_IS}), such that $A_\ev=0$, $C=0$. Since
$K\subset\fo\fo_{I\Pi}(n_{\ev}|n_{\od})$, it follows that
$K^{(1)}\subset \fo\fo^{(1)}_{I\Pi}(n_{\ev}|n_{\od})$. We also
have $K\simeq \fo_\Pi(n_\od)$, so, as it was shown in subsec.
\ref{fo_P-der}, $K^{(1)}$ consists of matrices from $K$, such that
$A_\od$ has the shape (\ref{**}) (even if $n_\od<6$). Hence,
$K^{(1)}, F^{ij}, G^{ij}, H^i$ span $M$, so $M\subset
\fo\fo_{I\Pi}^{(1)}(n_{\ev}|n_{\od})$, and $M=
\fo\fo_{I\Pi}^{(1)}(n_{\ev}|n_{\od})$. Since $F^{ij}, G^{ij},
H^{ij}, I_\od\in\fo\fo_{I\Pi}^{(1)}(n_{\ev}|n_{\od})$, and
$K\subset \fo\fo_{I\Pi}^{(1)}(n_{\ev}|n_{\od})$ (because $K^{(1)}$
and $H^i$ span $K$), we also get that $M\subset
\fo\fo_{I\Pi}^{(2)}(n_{\ev}|n_{\od})$, and $M=
\fo\fo_{I\Pi}^{(2)}(n_{\ev}|n_{\od})$. \end{proof}

\ssec{The derived Lie superalgebras of $\fo\fo_{\Pi\Pi}(n_{\ev}|
n_{\od})$} Now we consider the case where $n_\ev, n_\od$ are even
and $B$ is of the shape (\ref{canon_sup_ev}) such that
 $B_i=\Pi(n_i)$. We set $k_\ev=n_\ev/2$, $k_\od=n_\od/2$. We will denote the Lie superalgebra
preserving this form $B$ by $\fo\fo_{\Pi\Pi}(n_{\ev}| n_{\od})$;
it consists of supermatrices of the following shape:
\begin{equation}
\label{oo_SS}
\begin{pmatrix}
\Pi(n_\ev)A_{\overline{0}}&\Pi(n_\ev)C^T\Pi(n_\od)\\C&\Pi(n_\od)A_{\overline{1}}\end{pmatrix},
\;
\begin{matrix}\text{where
$A_{\overline{0}}, A_{\od}$ are symmetric, }\\
\text{$C$ is an arbitrary $n_{\overline{1}}\times
n_{\overline{0}}$ matrix.}\end{matrix}
\end{equation}

Direct computation shows that
$$
\fo\fo_{\Pi\Pi}^{(i)}(2|2)=\begin{cases} \{\text{matrices of the
shape~} (\ref{oo_SS}) \text{~such that~} A_\ev, A_\od\in ZD(2)\}&
\text{if~} i=1,\\
\{\text{matrices of the shape~} (\ref{oo_SS})
\text{~such that~} \Pi(2)A_\ev= \Pi(2)A_\od=\lambda\cdot 1_2\}&
\text{if~} i=2,\\
\{\lambda\cdot 1_{2|2}\mid \lambda\in\Kee\}&
\text{if~} i=3,\\
0 & \text{if~} i\geq 4.
\end{cases}
$$

\begin{Lemma}\label{oo_PP} If $n_\ev+n_\od\geq 6$, then

i) $\fo\fo_{\Pi\Pi}^{(1)}(n_{\ev}|n_{\od})$ consists of the
matrices of the shape (\ref{oo_SS}) such that $A_\ev, A_\od$ are
zero-diagonal;

ii) $\fo\fo_{\Pi\Pi}^{(2)}(n_{\ev}|n_{\od})$ consists of matrices
from $\fo\fo_{\Pi\Pi}^{(1)}(n_{\ev}|n_{\od})$ such that
$$
\sum\limits_{i=1}^{n_\ev/2}
(\Pi(n_\ev)A_\ev)_{ii}+\sum\limits_{i=1}^{n_\od/2}
(\Pi(n_\od)A_\od)_{ii}=0,
$$
i.e.,  the \lq\lq half-supertrace" of the matrix vanishes;

iii)
$\fo\fo_{\Pi\Pi}^{(3)}(n_{\ev}|n_{\od})=\fo\fo_{\Pi\Pi}^{(2)}(n_{\ev}|n_{\od})$.
\end{Lemma}

\begin{proof} Let $\tilde{B}=\diag_2(\Pi(n_\ev),\Pi(n_\od))$ be a (non-super)
bilinear form on a space of dimension $n_\ev+n_\od$. Denote:
$$
\renewcommand{\arraystretch}{1.4}
\begin{array}{l}
M^k=\{L\in\fgl(n_\ev|n_\od)\mid \text{~exists~}
L'\in\fo^{(k)}_{\tilde{B}}(n_\ev+n_\od) \text{~such that~}
L_{ij}=L'_{ij} \text{~for all~}
i,j\in\overline{1,n_\ev+n_\od}\};\\
 N=\text{Span}\{L^2\mid
L\in\fo\fo_{\Pi\Pi}(n_{\ev}| n_{\od}), L \text{~is
odd~}\}.\end{array}
$$

As it was noticed before, matrices and supermatrices behave
identically with respect to multiplication and Lie (super)bracket.
So we get the following inclusion:
$$
M^i\subset \fo\fo_{\Pi\Pi}^{(i)}(n_{\ev}|n_{\od})\subset M^i+N.
$$

Since $\tilde{B}$ is equivalent to $\Pi(n_\ev+n_\od)$, it follows
from Lemma \ref{o_p} that $M^1,M^2,M^3$ coincide with conjectural
spaces $\fo\fo_{\Pi\Pi}^{(1)}(n_{\ev}|n_{\od})$,
$\fo\fo_{\Pi\Pi}^{(2)}(n_{\ev}|n_{\od})$,
$\fo\fo_{\Pi\Pi}^{(3)}(n_{\ev}|n_{\od})$ as they are described in
the lemma. So, to prove the lemma, it suffices to show that
$N\subset M^2$. If
$$
L=\begin{pmatrix} 0&\Pi(n_\ev)C^T\Pi(n_\od)\\C&0\end{pmatrix}
$$
is an odd matrix from $\fo\fo_{\Pi\Pi}(n_{\ev}|n_{\od})$, then
$$
L^2=\begin{pmatrix}
\Pi(n_\ev)C^T\Pi(n_\od)C&0\\0&C\Pi(n_\ev)C^T\Pi(n_\od)\end{pmatrix},
$$
and
$$
(C^T\Pi(n_\od)C)_{ii}=\sum\limits_{j=1}^{k_\od}
C_{ji}C_{j+k_\od,i}+\sum\limits_{j=k_\od+1}^{n_\od}
C_{ji}C_{j-k_\od,i}=0;
$$
similarly, $\Pi(n_\od)C\Pi(n_\ev)C^T\Pi(n_\od)$ is zero-diagonal,
so $N\subset M^1$. Now,
$$
\begin{array}{ll}
\sum\limits_{i=1}^{k_\ev} (\Pi(n_\ev)C^T\Pi(n_\od)C)_{ii}+
\sum\limits_{l=1}^{k_\od} (C\Pi(n_\ev)C^T\Pi(n_\od))_{ll}=\\
\sum\limits_{i=1}^{k_\ev}(\sum\limits_{j=1}^{k_\od}
C_{ji}C_{j+k_\od,i}+\sum\limits_{j=k_\od+1}^{n_\od}
C_{ji}C_{j-k_\od,i})  +  \sum\limits_{l=1}^{k_\od}
(\sum\limits_{m=1}^{k_\ev} C_{lm}C_{l,m+k_\ev}+
\sum\limits_{m=k_\ev+1}^{n_\ev} C_{lm}C_{l,m-k_\ev})=0,
\end{array}
$$
so $N\subset M^2$. \end{proof}

\begin{Lemma} If $n_{\ev}+n_{\od}\geq 6$, then

i) if $k=(n_{\ev}+n_{\od})/2$ is odd, then
$\fo\fo_{\Pi\Pi}^{(2)}(n_{\ev}|n_{\od})$ is simple;

ii) if $k$ is even, then the only non-trivial ideal of
$\fo\fo_{\Pi\Pi}^{(2)}(n_{\ev}|n_{\od})$ is $Z=\{\lambda\cdot
1_{n_\ev|n_\od}\mid \lambda\in\Kee\}$, its center, and hence
$\fo\fo_{\Pi\Pi}^{(2)}(n_{\ev}|n_{\od})/Z$ is simple.\end{Lemma}

\begin{proof} Let $\iota: \fgl(n_\ev|n_\od)\to \fgl(n_\ev+n_\od)$ be a forgetful map
that sends a supermatrix into the matrix with the same entries and superstructure forgotten. Since
for $p=2$ matrices and supermatrices behave identically with respect to the Lie
(super)bracket, $\iota(\fo\fo_{\Pi\Pi}^{(2)}(n_{\ev}|n_{\od}))$ is
a Lie algebra, and for any ideal
$I\subset\fo\fo_{\Pi\Pi}^{(2)}(n_{\ev}|n_{\od})$, $\iota(I)$ is an
ideal in $\iota(\fo\fo_{\Pi\Pi}^{(2)}(n_{\ev}|n_{\od}))$.

Set
$$
M=\begin{pmatrix}1_{k_\ev}&0&0&0\\0&0&1_{k_\ev}&0\\0&1_{k_\od}&0&0\\
0&0&0&1_{k_\od}\end{pmatrix}.
$$

Then, according to Lemma \ref{o_p} and \ref{oo_PP}, the map
$X\mapsto MXM^{-1}$ gives us an isomorphism between
$\iota(\fo\fo_{\Pi\Pi}^{(2)}(n_{\ev}|n_{\od}))$ and
$\fo_\Pi(n_\ev+n_\od)$. So, according to Lemma \ref{o_p_sim}, if
$k$ is odd, then $\iota(\fo\fo_{\Pi\Pi}^{(2)}(n_{\ev}|n_{\od})$ is
simple; if $k$ is even, then the only non-trivial ideal of
$\iota(\fo\fo_{\Pi\Pi}^{(2)}(n_{\ev}|n_{\od}))$ is
$$
\{\lambda\cdot M^{-1}1_{n_\ev+n_\od}M=\lambda\cdot
1_{n_\ev+n_\od}\mid \lambda\in\Kee\}=\iota(Z).
$$

Thus, since $\iota$ is invertible,
$\fo\fo_{\Pi\Pi}^{(2)}(n_{\ev}|n_{\od})$ can have a non-trivial
ideal only if $k$ is even, and this ideal must be equal to $Z$;
direct computation shows that $Z$ is indeed an ideal.\end{proof}

\ssec{The derived Lie superalgebras of $\fpe(k)$} As it was shown in
subsec.~\ref{SCanon}, any non-degenerate odd symmetric bilinear
form on a superspace of dimension $(k|k)$ (if dimensions of the
even and odd parts of the space are not equal, such a form does
not exist) is equivalent to $\Pi_{2k}$.

The Lie superalgebra $\fpe(k)$ preserving this form consists of
the supermatrices of the shape
\begin{equation}
\label{pe1} \left(\begin{array}{ll} A&C\\D&A^T\end{array}\right),
\quad \text{where $A\in\fgl(k)$, $C$ and $D$ are symmetric
$k\times k$ matrices}. 
\end{equation}

Direct computations show that:
\medskip

$\fpe^{(1)}(1)=\{\lambda\cdot 1_2\mid \lambda\in\Kee\}$;

$\fpe^{(2)}(1)=0$;

$\fpe^{(1)}(2)=\{\text{matrices of the shape (\ref{pe1}) such that
$B$ and $C$ are zero-diagonal}\}$;

$\fpe^{(2)}(2)=\{\text{matrices of $\fpe^{(1)}(2)$ such that $\tr
A=0$}\}$;

$\fpe^{(3)}(2|2)=\{\lambda\cdot 1_{2}\mid \lambda\in\Kee\}$;

$\fpe^{(4)}(2)=0$.

\begin{Lemma} If $k\geq 3$, then

i) $\fpe^{(1)}(k)=\{\text{matrices of the shape (\ref{pe1}) such
that $B$ and $C$ are zero-diagonal}\}$;

ii) $\fpe^{(2)}(k)=\{\text{matrices of $\fpe^{(1)}(k)$ such that
$\tr A=0$}\}$;

iii) $\fpe^{(3)}(k)=\fpe^{(2)}(k)$.\end{Lemma}

\begin{proof} As it was noticed before, matrices and supermatrices
behave identically with respect to multiplication and Lie
(super)bracket. So, if we denote
$$
\renewcommand{\arraystretch}{1.4}
\begin{array}{l}
M^i=\{L\in\fgl(k|k)\mid \text{~exists~} L'\in\fo_\Pi(2k)
\text{~such that~} L_{ij}=L'_{ij} \text{~for all~}
i,j\in\overline{1,2k}\};\\
 N^i=\text{Span}\{L^2\mid
L\in\fpe^{(i-1)}(k), L \text{~is odd~}\},\quad  \text{where~}
\fpe^{(0)}(k)=\fpe(k),\end{array}
$$
then we get the following inclusion:
$$
M^i\subset \fpe^{(i)}(k)\subset M^i+N^i.
$$

Now recall from the proof of Lemma \ref{o_p} that
$M^1,M^2,M^3$ coincide with conjectural spaces $\fpe^{(1)}(k)$,
$\fpe^{(2)}(k)$, $\fpe^{(3)}(k)$ as they are described in the
lemma. Notice also that $N^1\subset M^1$ (since $N^1$ is an even
subspace of $\fpe(k)$), so $\fpe^{(1)}(k)=M^1$. If
$$
L=\left(\begin{array}{ll}0&C\\D&0\end{array}\right)
$$
is an odd matrix from $\fpe^{(1)}(k)$ (so $C$ and $D$ are
symmetric and zero-diagonal), then
$$
L^2=\left(\begin{array}{ll}CD&0\\0&DC\end{array}\right),
$$
and
$$
\tr CD=\sum\limits_{i,j=1}^k C_{ij}D_{ji}=2\sum\limits_{1\leq
i<j\leq n} C_{ij}D_{ij}=0,
$$
so $N^2\subset M^2$, and $\fpe^{(2)}(k)=M^2$. Also, $N^3\subset
N^2\subset M^2=M^3$, and $\fpe^{(3)}(k)=M^3$. \end{proof}

\section{Canonical expressions of non-symmetric bilinear superforms. Related
Lie superalgebras.} As it was shown in sec.
\ref{non-sym_super}, any even non-symmetric bilinear form $B$ on a
superspace in the
 standard format has the shape
$$
\left(\begin{array}{ll} B_\ev&0\\
0&B_\od\end{array}\right).
$$
This form is preserved by the matrices of the shape
$$
\left(\begin{array}{ll} A_\ev&C\\
D&A_\od\end{array}\right),\text{~where~}\begin{cases}
A_\ev \text{~preserves the form~} B_\ev,&\\
A_\od \text{~preserves the form~} B_\od,&\\
CB_\od+B_\ev D^T=0,&\\
DB_\ev+B_\od C^T=0.\end{cases}
$$

Any odd non-degenerate non-symmetric bilinear form $B$ on a
superspace of dimension $(k|k)$ (non-degenerate odd form can exist
only on a superspace with equal even and odd dimensions) is
equivalent to a form of the shape
$$
\left(\begin{array}{ll} 0&1_k\\
J&0\end{array}\right),\text{~where~} J \text{~is the Jordan normal
form}.
$$
Such a form is preserved by the matrices of the shape
$$
\left(\begin{array}{ll} A&C\\
D&A^T\end{array}\right),\text{~where~}\begin{cases}
 AJ^T+J^TA=0,&\\
CJ+C^T=0,&\\
D+JD^T=0.\end{cases}
$$

\begin{Problem} The (first or second) derived Lie superalgebra of $\fpe_B(k|k)$
is simple (perhaps, modulo center) only if $J$ consists of $1\times 1$ blocks. What is an explicit structure
in other cases? Compare with Ermolaev's description \cite{Er}. \end{Problem}

\end{document}